  \documentclass{elsart}

  \usepackage{graphics}
  \usepackage{amssymb}
\usepackage{epsfig}

\def\JCPfixesON{}
\newcommand{\JCPvspace}[1]{%
 \ifx\JCPfixesON\undefined 
 \else \vspace*{#1} 
 \fi}
\newcommand{\JCPnewpage}{%
 \ifx\JCPfixesON\undefined 
 \else \newpage 
 \fi}

\newlength{\RRRfigCapWid}
\ifx \epsfig\undefined \newcommand{\epsfig}[1]{} \fi
\ifx  \psfig\undefined \newcommand{\psfig}[1]{}  \fi
\newcommand{\figPIC}[3]{
 \ifx\USEepsfig\undefined
    \ifx\USEpsfig\undefined
       \fbox{\parbox{#2}{\rule{0in}{#3}}}
    \else
       \mbox{\psfig{figure=#1,height=#3,clip=}}
    \fi
 \else
    \epsfig{figure=#1,height=#3}
 \fi}

\newcommand{\makemylabel}[2]{
   \ifx\theCOUNTERusedBYmakemylabel\undefined 
      \newcounter{COUNTERusedBYmakemylabel}   \fi 
   \renewcommand{\theCOUNTERusedBYmakemylabel}{#2}
   \refstepcounter{COUNTERusedBYmakemylabel}   \label{#1}}
\newcommand{\GetRemLength}[1]{
   \ifx\AuxRemLength\undefined  \newlength{\AuxRemLength}  \fi
   \ifx\RemLength\undefined     \newlength{\RemLength}     \fi
   \settowidth{\AuxRemLength}{#1}
   \setlength{\RemLength}{\textwidth}
   \addtolength{\RemLength}{-\AuxRemLength}}

\newcommand{\dvn}[3]{\frac{d^{#3} #1}{d #2^{#3}}}

\newcommand{\eps}{\varepsilon}

\newcommand{\bfrac}[2]{\left(\frac{#1}{#2}\right)}

\newcommand{\RRRfigOLD}[3]{\begin{figure} [htb!] \begin{center}
   \ifx\USEepfig\undefined #1 \fi  \ifx\USEpsfig\undefined #2 \fi  #3
   \end{center} \end{figure}}

\newcommand{\EPfg}[1]{\epsfig{figure=./#1}}

\newcommand{\Dmath}[1]{${\displaystyle #1}\/$}

\newcommand{\dst}[1]{{\displaystyle #1}}
\newtheorem{PWFalgorithm}{Algorithm}[section]

\usepackage{amsmath}
\usepackage{enumerate}

\begin{document}
  \begin{frontmatter}

\title{Multirate integration of axisymmetric step-flow equations}
\author{Pak-Wing Fok\corauthref{cor}}
\address{Applied and Computational Mathematics,
         California Institute of Technology,
         Pasadena, CA 91125}
\address{Department of Biomathematics,
         University of California, Los Angeles,
         CA 90095}
\corauth[cor]{Corresponding author.}
\ead{pakwing@caltech.edu}
\author{Rodolfo R. Rosales}
\address{Department of Mathematics,
         Massachusetts Institute of Technology,
         Cambridge, MA 02139}

\begin{abstract}
 \noindent
 We present a multirate method that is particularly suited for
 integrating the systems of Ordinary Differential Equations (ODEs)
 that arise in step models of surface
 evolution. The surface of a crystal
 lattice, that is slightly miscut from a plane of symmetry,
 consists  of a series of terraces separated by steps. Under the
 assumption of axisymmetry, the step radii satisfy a system of ODEs
 that reflects the steps' response to step line tension and
 step-step interactions. Two main problems arise in the numerical
 solution of these equations. First, the trajectory
 of the innermost step can become singular, resulting in
 a divergent step velocity. Second, when
 a step bunching instability arises, the motion of
 steps within a bunch becomes very strongly stable, resulting in
 ``local stiffness''. The multirate method introduced in this paper
 ensures that small time steps are taken for
 singular and locally stiff components, while larger time 
 steps are taken for the remaining ones. 
 Special consideration is given to the construction
 of high order interpolants during run time which ensures
 fourth order accuracy of scheme for components of the
 solution sufficiently far away from singular trajectories.
\end{abstract}

  \end{frontmatter}

\section[
Introduction%
]{
Introduction.%
}\label{sec01:introduction}
In this paper, we present a multirate method that is suited for
integrating systems of Ordinary Differential Equations (ODEs)
which arise from step models describing nanostructure evolution.
The strength of our method, in comparison to other existing
multirate schemes, is its order of accuracy: our method
is fourth order provided solutions are sufficiently smooth. 
%
%
A multirate method is
basically one that takes different step sizes
for different components of the solution \cite{GearWellsBIT84}. When might
such a need for different time steps arise? One situation where
multirate methods may be more efficient
than single rate ones is when a few of the components
contain time singularities or are \emph{locally stiff}. In this case,
(explicit) single rate methods are likely to use small time steps
for all the components, whereas a multirate one uses small
time steps just for singular/locally stiff ones. 
We believe that this strategy 
significantly improves the efficiency of solution.
Integration by our multirate method occurs in two stages.
We first use an ``outer'' integrator to handle the non-singular/
non-stiff components
and then an ``inner'' integrator to handle the singular/stiff ones. To perform
the second inner integration, a small number of components from
the outer solution must be interpolated. The high order of accuracy
of the multirate scheme relies on the interpolation being of
a sufficiently high order. One of our key results in this paper
is that for a multirate method to be $n^{th}$
order, the interpolation must be $(n-1)^{st}$ order. We demonstrate
our method by using a fourth order Runge-Kutta method 
(with error control) for the inner and outer
schemes and coupling them together with cubic interpolants. 

Multirate schemes were first studied by Gear and
Wells~\cite{GearWellsBIT84}. Further treatments can be found
in~\cite{LoggApplNumMath04,SavcencoHundsdorferBIT06,MakinoAarsethPublAstronSocJap92,WaltzPageJCP02,KatoKataokaElecEngJap99,BartelGuntherJCompApplMath02,GuntherKvaernoBIT41}, and the
references therein. However, considering the wide variety of methods
which researchers have used to improve the performance and accuracy
of integration codes, it is surprising to learn that multirate
schemes have received only modest attention.
This is even more surprising, given that algorithms using somewhat
similar concepts are well developed for the numerical solution of
PDEs --- e.g. Adaptive Mesh Refinement (AMR) in the field of
Hyperbolic PDEs~\cite{BergerOligerJCP84,HenshawSchwendemanJCP03}.
The applications of multirate methods seem mostly confined to
$N$-body problems~\cite{MakinoAarsethPublAstronSocJap92,%
WaltzPageJCP02}, and equations arising from electrical
networks~\cite{KatoKataokaElecEngJap99,BartelGuntherJCompApplMath02}.
The work presented here is, as far as we are aware, the first
instance of a multirate scheme applied to a problem in surface
evolution.

Our method is fairly similar to the one described
in~\cite{SavcencoHundsdorferBIT06}, which is second order.
For example: we first advance and
interpolate the slow components, and then integrate the fast
components --- this is the ``slowest first'' paradigm described
in~\cite{GearWellsBIT84}. We also automatically detect the fast
components by looking at the error estimates produced using an
embedded formula\footnote{Some methods rely on the user knowing
enough about the physical system at hand, so that the fast and
slow components are known in advance \cite{GuntherKvaernoBIT41}.}.
The main differences are that: (i) our method is
fourth order in time for sufficiently smooth solutions and 
(ii) our method is most effective when applied to the types
of ODEs that commonly arise in models of step surface evolution.
We will discuss (ii) in more detail later on in the paper, but
for now we make the comment that a necessary but not sufficient
condition for our method to work is that the ODE system be
locally coupled. The strengths of our method are that automatic step size
selection is simple to implement, and that there is a lower
overhead cost because we do not have to interpolate all the slow
components of the solution --- this is one of the benefits of
specializing to locally coupled systems.
%
%

The outline of this paper is as follows: 
in Section \ref{sec:physical_motivation}, we 
discuss why step models are studied and
explain the physics behind the equations.
In Section
\ref{sec02:governing} we present the step equations, and in
Section \ref{sec03:properties} we describe features of the 
equations that require special attention. In Section \ref{sec04:code}
we give the details of the numerical method and in particular,
we discuss the important
issue of interpolation in Section \ref{interpolation}. 
We validate our code and present our results in Section 
\ref{sec05:implementation} and summarize our findings with a conclusion
in Section \ref{sec61:discussion}.

\section[
Physical Motivation%
]{
Physical Motivation.%
}\label{sec:physical_motivation}
The surface dynamics of crystal structures has received much recent
attention~\cite{JeongWilliamsSurfSciRep99,ThurmerReutt-RobeyPRL01,%
YagiMinodaSurfSciRep01} because of its relevance to the
fabrication of nano-scale electronic devices, such as quantum
dots~\cite{KitamuraNishiokaApplPhysLett95}. Of interest to us here
is the behavior of \emph{vicinal} surfaces below the
roughening temperature $T_R\/$. For any given material, a
surface forming a small angle with a high 
symmetry plane of the crystal is called
vicinal. The roughening temperature is the critical temperature below
which steps become thermodynamically stable. Thus, a microscope image
\RRRfigOLD{\EPfg{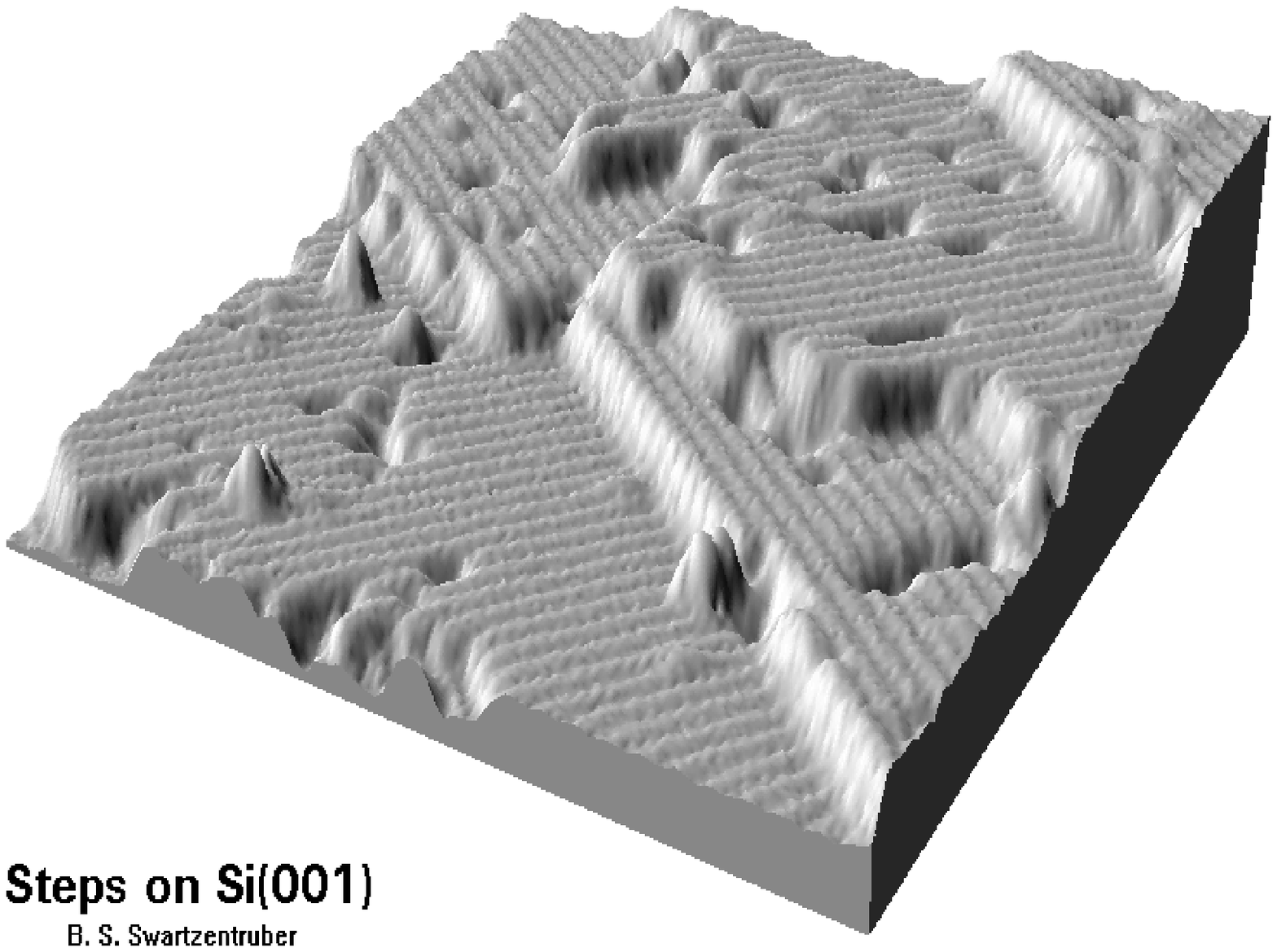,width=4in}}{
\caption{Image of a step on Silicon (001), taken with a Scanning
Tunneling Microscope (courtesy of Brian Swartzentruber, 
Sandia National Laboratory).} 
\label{fig:silicon}}

of a surface below $T_R$, with a slight miscut angle, appears as made
from a series of terraces separated by steps of atomic height --- see
Figure~\ref{fig:Silicon} for an example in Silicon.
As the surface evolves, the steps move and change their shape, but
the steps are well defined and have a lifetime that is long enough to
be directly observable. When the temperature is increased above
$T_R$, a \emph{Kosterlitz-Thouless} phase transition
occurs~\cite{ChuiWeeksPRB01}, and the surface becomes statistically
rough --- as characterized by the divergence of the height--height
correlation
function~\cite{PrasadWeichmanPRB98,VillainGrempelJPhysF85}. For many
physical applications (such as epitaxy) the operating temperatures
are below $T_R\/$, and this mesoscale description of a surface in
terms of steps and their evolution is very useful. A step model can
account for finite size effects occurring at the atomic scale, while
remaining computationally simple --- simulations with step models
can be done over much longer time periods than, for example, with
atomistic models of the surfaces.

The step's shape can, in general, be very complicated. Thus,
quantitative descriptions for how steps interact with one another
can be very difficult to derive, and a complete description of an
arbitrarily shaped nanostructure in terms of its steps is currently
not available. As a result, theoretical studies have been restricted
to simple nanostructure geometries and step shapes. The BCF model,
proposed by Burton, Cabrera and
Frank~\cite{BurtonCabreraPhilosTrans51} in 1951, deals with a
monotonic step train consisting of an infinite number of parallel
steps --- see Figure~\ref{fig:GeomBCF}.
\RRRfigOLD{\EPfg{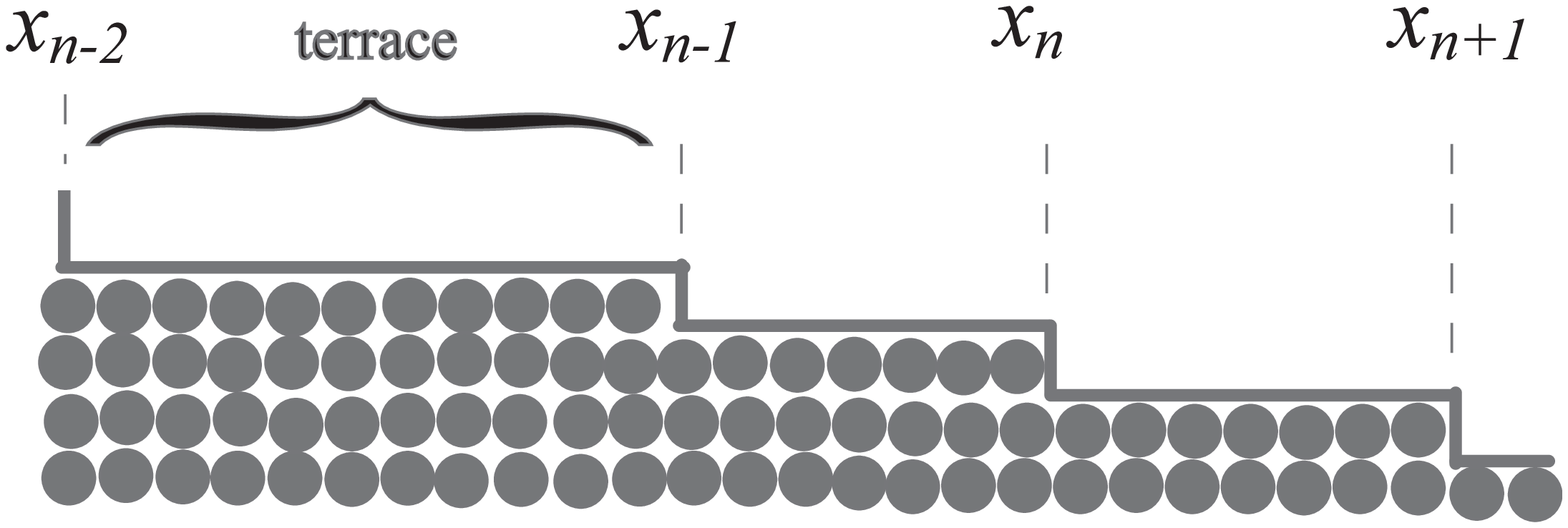,width=3.5in}}{
 \caption{Sketch of the geometry in the Burton, Cabrera and Frank (BCF) 1-D step
model. The step positions $x_n\/$ and the terraces are shown.}
\label{fig:GeomBCF}}
%
%
In this model, the steps edges are separated by atomically smooth
terraces, and each step position is uniquely described by a single
scalar quantity $x_n(t)\/$ --- where the index $n\/$ labels the step.
Using this model, Burton, Cabrera and Frank were able to describe
the evolution of a (1 dimensional) stepped surface, under
non-equilibrium conditions, in terms of its steps.

In 1988, Rettori and Villain \cite{RettoriVillainJPhys88} considered
a 2D array of circular mounds, and incorporated the effects of step
line tension into the BCF model. The nanostructures that they studied
consist of a finite number of concentric circular layers, in a
``wedding cake'' configuration -- see Figure~\ref{fig:WedCake}. 
This step system can also be used to describe the ``healing''
of small circular pits \cite{YamamotoSudohSurfSci06} produced 
by scanning tunneling microscopes. The
\RRRfigOLD{\EPfg{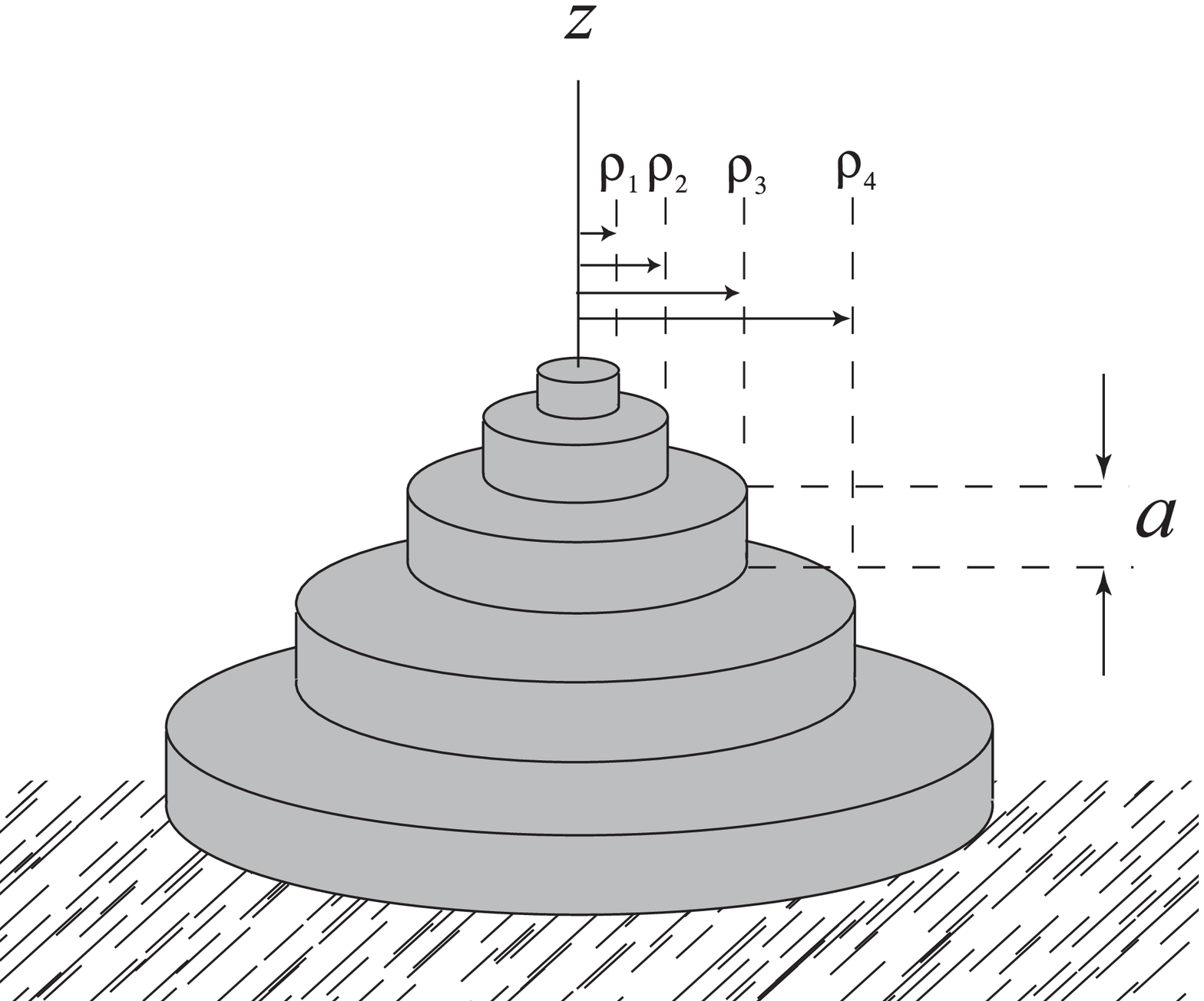,width=3in}}{
\caption{
``Wedding Cake'' step configuration for an axisymmetric nanostructure,
with a finite number of concentric, circular steps. The step
height $a$ is of the order of the crystal's lattice constant.
The number of steps considered in this paper is much larger
than what is shown in this picture.} 
\label{fig:WedCake}}

radius of each layer $\rho_n(t)\/$ is assumed to be a continuous
function of time. Physical considerations then lead to a set of
locally coupled ODEs for the radii. Similar sets
of equations can be found in 
\cite{IsraeliKandelPRB99,MargetisAzizPRB05,SatoUwahaSurfSci99}.
There are two main competing physical processes that take place on a
stepped surface, in the absence of evaporation and desorption. The
first one is the diffusion of adsorbed atoms (``adatoms'') across
terraces, which is characterized by a diffusivity $D_s\/$. The second
one is the attachment-detachment of adatoms at step edges, which is
characterized by the kinetic coefficients $k_{+}\/$ and $k_{-}\/$ ---
see Figure~\ref{fig:TranProStepSurf}.
Experimental evidence~\cite{EhrlichHuddaJChemPhys66} suggests that,
for some materials, attachment from the terrace above requires
overcoming a higher activation energy barrier than attachment from
the terrace below, so that $k_{+} > k_{-}\/$. However, in this paper
we consider only $k_{+} = k_{-} = k\/$, and disregard this (possible) asymmetry
in the step attachment-detachment --- known in the literature as an
Ehrlich-Schwoebel (ES) Barrier~\cite{EhrlichHuddaJChemPhys66,%
SchwoebelJApplPhys69,SchwoebelShipseyJApplPhys66}. Furthermore, we
neglect the diffusion of adsorbed vacancies and the diffusion
of adatoms along step edges \cite{MargetisPRB07}.
\RRRfigOLD{\EPfg{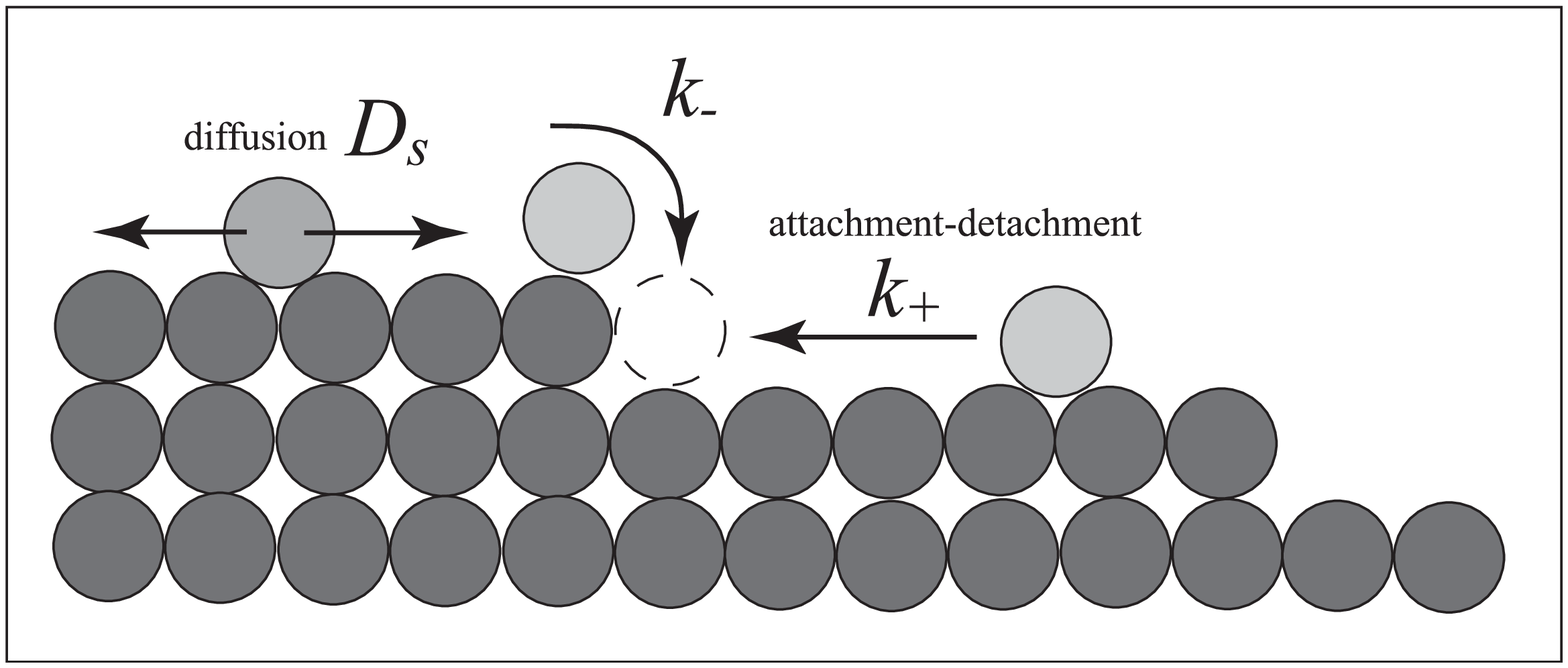,width=4in}}{
\caption{Sketch of the transport processes on a stepped surface, characterized
by the attachment-detachment coefficients $k_{+}\/$ and $k_{-}\/$, and
by the terrace diffusivity $D_s\/$.
Throughout this paper, we take $k_{+} = k_{-}$.} 
\label{fig:TranProStepSurf}}


\section[
Governing Equations%
]{
Governing Equations.%
}\label{sec02:governing}
In this section, we present the (non-dimensional) evolution equations
for a finite axisymmetric nanostructure with \Dmath{N \geq 5\,}
steps, relaxing in the absence of deposition and evaporation
--- see Figure~\ref{fig:WedCake}. Derivations of these equations can
be found in \cite{IsraeliKandelPRB99,MargetisAzizPRB05}.

Every step in the structure is subject to two physical effects that
drive its motion. The first is a \emph{step-line tension}, arising
from a Gibbs-Thomson mechanism~\cite{McleanKrishnamachariPRB97}. An
isolated, circular step of radius $\rho(t)\/$, on top of an infinite
substrate, initially devoid of adatoms, reduces its perimeter (and hence
its radius) by emitting adatoms at a rate proportional to its
curvature~\cite{JeongWilliamsSurfSciRep99} --- i.e.
\Dmath{\dot{\rho} \propto -1/\rho}. The second effect is a
\emph{repulsive interaction with neighboring steps},
characterized by a potential function that is inversely proportional
to the square of the distance between the steps~\cite{TanakaBarteltPRL97}.
Steps in the bulk of the structure (with a smaller curvature) tend to
be less affected by the step-line tension compared to steps near the top. 

Let \Dmath{\rho_n = \rho_n(t)} \Dmath{(1 \leq n \leq N)} be the
positions\footnote{Non-dimensionalized radii, measured from the axis
   of symmetry.}
of the steps --- numbered starting from the top of the
nanostructure. Thus
\Dmath{0< \rho_1 < \rho_2 \dots < \rho_N < \infty}, with
\Dmath{\rho_1} the radius of the innermost step and \Dmath{\rho_N}
the radius of the outermost step. Define \Dmath{\Lambda_n},
\Dmath{R_n}, and \Dmath{\Delta_n}, by:
\begin{equation} \label{eqn:sec20:Lambdan}
  \parbox[c]{4.6in}{\Dmath{\left. \begin{array}{lclcl}
     \Lambda_1 & = & \lambda(\rho_{1}, \rho_{2})\/,\\
     \Lambda_n & = & \lambda(\rho_{n-1}, \rho_{n}) +
                     \lambda(\rho_{n}, \rho_{n+1})\/,
                     \quad \mbox{for} \;\; 2 \leq n \leq N-1\/,\\
     \Lambda_N & = & \lambda(\rho_{N-1}, \rho_{N})\/,
  \end{array} \; \right\}}}
\end{equation}
\begin{equation} \label{eqn:sec20:lambda}
  \parbox[c]{4.6in}{\Dmath{
     \hspace*{0.50in} \mbox{where} \quad
     \lambda(\rho_i, \rho_j) = \frac{2\,\rho_i}{\rho_i+\rho_j}\,
     \frac{1}{(\rho_i-\rho_j)^3} + \frac{1}{\rho_j}\,
     \left(\frac{\rho_i}{\rho_i^2-\rho_j^2}\right)^2\,,}}
\end{equation}
\begin{equation} \label{eqn:sec20:Rn}
  \parbox{4.4in}{\Dmath{\left. \begin{array}{lcl}
     R_n & = & \dst{\frac{1}{\rho_n}} + \varepsilon\,\Lambda_n\/,
     \quad \mbox{for} \;\; 1 \leq n \leq N\/,
  \end{array} \right.}}
\end{equation}
\begin{equation} \label{eqn:sec20:Deltan}
  \parbox{4.4in}{\Dmath{\left. \begin{array}{lcl}
     \Delta_n & = &
          m_1\,\ln \dst{\frac{\rho_{n+1}}{\rho_n}} + m_2\,\left(
          \dst{\frac{1}{\rho_{n+1}} + \frac{1}{\rho_{n}}} \right)\/,
     \quad \mbox{for} \;\; 1 \leq n \leq N-1\/,
  \end{array} \right.}}
\end{equation}
where \Dmath{\varepsilon}, \Dmath{m_1}, and \Dmath{m_2} (as well as
\Dmath{\gamma} below) are non-dimensional constants. Then the
step-flow equations are
\begin{eqnarray}
 \frac{d}{dt}\rho_1 & = & \frac{\gamma}{\rho_1}\,\left(
 \frac{R_{2}-R_{1}}{\Delta_{1}} \right)\/,
 \label{eqn:sec20:eqnrho1} \\
 \frac{d}{dt}\rho_n & = & \frac{\gamma}{\rho_n}\,\left(
 \frac{R_{n+1}-R_{n}}{\Delta_{n}} - \frac{R_{n}-R_{n-1}}{\Delta_{n-1}}
 \right)\/, \quad \mbox{for} \;\; 2 \leq n \leq N-1\/,
 \label{eqn:sec20:eqnrhon} \\
 \frac{d}{dt}\rho_N & = & \frac{\gamma}{\rho_N}\,\left(
  - \frac{R_{N}-R_{N-1}}{\Delta_{N-1}} \right)\/.
 \label{eqn:sec20:eqnrhoN}
\end{eqnarray}
The non-dimensional constants are as follows:
\begin{itemize}
  \item \vspace*{-1mm} 
  The parameter \Dmath{\eps > 0} measures the strength of the
  step-step interactions relative to the strength of the step line
  tension. It is given by
  \[
     \eps = \frac{2}{3}\,\frac{g_3}{g_1}\,\bfrac{a}{L}^2\/,
  \]
  where \Dmath{g_3} is the
  step-step interaction coefficient~\cite{TanakaBarteltPRL97},
  \Dmath{g_1} is the step stiffness~\cite{JeongWilliamsSurfSciRep99},
  \Dmath{a} is the height of a single step, and \Dmath{L} is a typical
  value for the radii -- for example, it could be the initial
  radius of the final step in the structure.
  \emph{We note that, in many experimental situations, \Dmath{0 < \eps \ll 1}}.
  \item 
  The parameters \Dmath{0 \leq m_1, m_2 \leq 1} are given by
  \[
     m_1 = \frac{k\,L}{k\,L+D_s} \quad \mbox{and} \quad
     m_2 = \frac{D_s}{k\,L+D_s}\/,
  \]
  where \Dmath{k} is the adatom attachment-detachment coefficient at
  a step, and \Dmath{D_s} is the adatom terrace diffusivity. The
  ratio
  \[
    m = \frac{m_2}{m_1\,\Delta \rho} = \frac{D_s}{k\,L_w}\/,
  \]
  where \Dmath{\Delta \rho = \rho_{n+1} - \rho_n\,} and \Dmath{L_w \equiv L \Delta \rho}
  is a typical terrace width, measures the competition
  between diffusion and attachment-detachment --- see equation
  (\ref{eqn:sec20:Deltan}).
  \item 
  The dimensionless parameter \Dmath{\gamma} is given by
  \[
     \gamma = \bfrac{g_1 \Omega_s}{k_B T} (\Omega_s c_s)
              \bfrac{a}{L}\bfrac{D_s}{LU}\,m_1\/,
  \]
  where \Dmath{\Omega_s} is the atomic area, \Dmath{k_B} is the
  Boltzmann constant, \Dmath{T} is the absolute temperature, \Dmath{c_s} is
  the equilibrium density of adatoms at a straight, isolated, step
  and \Dmath{U} is a typical bulk step velocity.
\end{itemize}
We note that, of the physical parameters involved in the definitions
of the non-dimensional constants above, some --- such as the terrace
diffusivity \Dmath{D_s}, have been extensively tabulated
\cite{KelloggSurfSciRep94}, while others --- such as \Dmath{c_s} and
\Dmath{U}, can be inferred from experiments
\cite{ThurmerReutt-RobeyPRL01,IchimiyaTanakaPRL96}. However, for the
purposes of simulation, we can take
\Dmath{m_1 = m_2 = \gamma = 1/2} without loss of generality, by an
appropriate rescaling of the step radii and time.

Equations (\ref{eqn:sec20:eqnrho1}) -- (\ref{eqn:sec20:eqnrhoN})
constitute a pentadiagonal system, of the form
\begin{equation} \label{eqn:sec20:pentadiagonal}
   \frac{d\rho_n}{dt}  = F_n \left(  \rho_{n-2}\/,\,\rho_{n-1}\/,\,
   \rho_n\/,\,\rho_{n+1}\/,\,\rho_{n+2} \right)\/,
\end{equation}
where all the \Dmath{F_n} have the same functional form,
with the exception of \Dmath{F_1}, \Dmath{F_2}, \Dmath{F_{N-1}}, and
\Dmath{F_N} --- which govern the behavior of the first and final two
steps.
Notice that if \Dmath{F = F_n} denotes the common rate function for
the bulk of the steps (\Dmath{2<n< N-2}), then
\begin{eqnarray*}
  F_{N-1}(\rho_{N-3}\/,\,\rho_{N-2}\/,\,\rho_{N-1}\/,\,\rho_{N})
  & = & \lim_{\zeta \to \infty}
  F(\rho_{N-3}\/,\,\rho_{N-2}\/,\,\rho_{N-1}\/,\,\rho_{N}\/,\,\zeta)\/,
  \\
  F_{N}(\rho_{N-2}\/,\,\rho_{N-1}\/,\,\rho_{N})
  & = & \lim_{\eta \to \infty} \lim_{\zeta \to \infty} 
  F(\rho_{N-2}\/,\,\rho_{N-1}\/,\,\rho_{N}\/,\,\eta\/,\,\zeta)\/.
\end{eqnarray*}
In this paper, we present a multirate method for integrating the equations
(\ref{eqn:sec20:eqnrho1}) -- (\ref{eqn:sec20:eqnrhoN}), when
\Dmath{\eps \ll 1}. However, we believe that the method presented here should
be applicable to \emph{general sets} of locally coupled ODEs, which are
\emph{locally stiff} (see Section \ref{local}). 
In fact, our multirate method was designed to specifically
tackle this problem.

\section[
Properties of the step-flow equations%
]{
Properties of The Step-Flow Equations.
}\label{sec03:properties}
%

Let us now turn our attention to the difficulties that arise
when solving equations (\ref{eqn:sec20:eqnrho1}) -- (\ref{eqn:sec20:eqnrhoN})
numerically. The axisymmetric step-flow 
equations possess a number of peculiar properties which 
pose problems for standard
integrators -- hence the need
for a multirate algorithm. For example,
the singular collapse of the innermost step 
and stiffness, localized to only a few of
the components, are two (different)
situations under which a standard integrator is forced
to use small time steps. In these cases, a single
rate method uses small time steps for all components.
In contrast, a multirate method will use small time steps
only when it has to, so that most of the components are
integrated with a large time step. This strategy improves the
efficiency of the integration.

The singular collapse of the innermost step causes
a loss in accuracy for most high order integrators
near the point of collapse.
Hence, we implement a low-order Simple Euler routine 
for the innermost step and its neighbors. Away from the
singular step, we implement a fourth order multirate
method with error control, which is able to efficiently
integrate locally stiff components.
\subsection[
Singular Collapse of Steps%
]{
Singular Collapse of Steps.%
}\label{singular}
Equations (\ref{eqn:sec20:eqnrho1}) -- (\ref{eqn:sec20:eqnrhoN}) have
the property that \Dmath{\rho_1 \rightarrow 0} in a finite time. The top
step always undergoes a monotonic collapse because 
its radius always decreases under the effect of step-line tension.
As the top step shrinks, it
emits adatoms, causing the radii of the second and subsequent steps to
grow as these are absorbed. When the top step completely disappears,
the number of layers in the structure is reduced by one. As a result
of the sequential collapse of top steps, a macroscopically flat region
called a facet forms and grows on the top of the structure. Provided
that the collapse of the top steps is tracked accurately, and the
topmost \Dmath{\rho_i} is removed at each collapse, the growth of the
facet is automatically accounted for.

When the first collapse occurs, $\rho_2$, $R_3$, $R_2$ and $\Delta_2$
replace $\rho_1$, $R_2$, $R_1$ and $\Delta_1$ in (\ref{eqn:sec20:eqnrho1}),
and (\ref{eqn:sec20:eqnrhon}) applies when $3 \leq n \leq N-1$. A similar
replacement occurs for the second and subsequent collapses.
In this way, a given index tracks always the same step  throughout
the integration.

Let \Dmath{\tau_1} be the collapse time for \Dmath{\rho_1}. Then it can
be shown~\cite{MargetisPrivateCommunication} that as $t \rightarrow \tau_1$,
\begin{equation} \label{square_root}
   \rho_1 \sim C_1 (\tau_1 - t)^{1/2} +
               C_2 (\tau_1 - t) \, \ln (\tau_1 - t) + O(\tau_1-t)\/.
\end{equation}
For some constants $C_1$ and $C_2$. The square root behavior
in (\ref{square_root}) comes from the fact that the leading
order behavior for $\rho_1 \ll 1$ 
in equation (\ref{eqn:sec20:eqnrho1})
stems from a line tension: $\dot{\rho}_1 \sim 1/\rho_1$.
Thus the derivatives of \Dmath{\rho_1} are divergent at the time of
collapse. Since (\ref{eqn:sec20:eqnrho1}) -- (\ref{eqn:sec20:eqnrhoN})
is a locally coupled set of equations, we also expect
$\rho_i$, for $2 < i\ll N$, to be singular but for the solutions to become
more regular near $\tau_1$ as $i$ becomes larger.
Since the accuracy of high order integrators usually relies
on the solution having enough bounded derivatives, this means that
standard high order solvers will lose accuracy
near the time of collapse. For example, consider a method
with truncation error \Dmath{O((\Delta t)^p \, y^{(p)})} for smooth
solutions \Dmath{y = y(t)}, and time step \Dmath{\Delta t}. Let
\Dmath{t} be a time \Dmath{m} steps away from
$\tau_1$, so that \Dmath{\tau_1-t = m\,\Delta t}. Then, given
the square root singularity in
(\ref{square_root}), the error near $\tau_1$ will be increased to
\begin{equation} \label{eq100}
   \mbox{error} = O(m^{-p+1/2} (\Delta t)^{1/2})\/.
\end{equation}
Furthermore, consider the issue of automatic time step selection in
an adaptive integration code. This is usually done by estimating the
local truncation error, and updating the time step size with a formula
for the truncation error that assumes a smooth solution. For example,
consider a Runge-Kutta scheme using an embedded higher order formula
to estimate the local truncation error. Such an algorithm updates the
time step size using a formula like \cite{Press}
\begin{equation} \label{stepsizecontrol}
   \Delta t_{\textrm{new}} = \Delta t_{\textrm{old}} \, \left( \frac{
   \textrm{desired error}}{\textrm{estimated error}}\right)^
   {\displaystyle \frac{1}{1 + \alpha}}\/,
\end{equation}
where \Dmath{\alpha} is the order of the integrator. Equation
(\ref{stepsizecontrol}) is invalid
near a time singularity, because from 
(\ref{eq100}), the error near $\tau_1$ does not
scale as \Dmath{\Delta t^{\alpha+1}}. The resulting behavior is
somewhat unpredictable: an adaptive integrator
may take a very large number 
of tiny steps -- rendering it very inefficient -- or it may
simply abort, stating that the specified error tolerance is not
achievable. 

\subsection[
Local Stiffness%
]{
Local Stiffness%
}\label{local}
The aim of this subsection is to attempt to quantify the classes of
systems for which the approach in this paper is effective.
Before explaining what we mean by local stiffness, we first
introduce some notation. For $i=1,2,...,N$, let $\rho_i(t)$
be the solution of the ODE system for some initial 
condition $\rho_i(0)$. For some integer $p$, let $\tilde{\rho}_i(t)$ 
be the solution with initial condition
$\rho_i(0)+\nu \delta_{p,i}$, 
where \Dmath{\delta_{p\/,i}} is the Kronecker delta, and \Dmath{\nu} is
small. Finally, let
$\rho(t) \equiv (\rho_1(t),\rho_2(t),...,\rho_N(t))$ and
$\tilde{\rho}(t) \equiv (\tilde{\rho}_1(t),\tilde{\rho}_2(t), ... ,\tilde{\rho}_N(t))$.
We say that 
the $p^{th}$ component of the ODE system is \emph{strongly local}
if (i) the system is locally coupled 
and (ii) given any $\epsilon > 0$, for all $t$, there exists an integer $d$ independent
of $t$ such that for $|n-p| > d$, $|\tilde{\rho}_n(t) 
- \rho_n(t)| < \epsilon$. If every component is strongly local, we say that
the ODE \emph{system} is strongly local.
Therefore, a system is strongly local if a small perturbation to any one of its components
remains localized in component number.

Now we explain local stiffness.
Recall that an ODE is stiff when 
the ratio of the slowest and fastest
time scales is much greater than one.
The simplest example of this is
a situation where the solution of interest is strongly stable, so that
small perturbations decay very rapidly, relative to the principal
time scale of evolution.
Now consider again the perturbation described in the previous paragraph.
We say that the $p^{th}$ component of an ODE system 
is \emph{locally stiff} if (i) it is strongly
local and (ii) $||\tilde{\rho}(t) - \rho(t)|| \rightarrow 0$ 
rapidly in time, relative to the principal
time scale of evolution. 
Hence, the $p^{th}$ component of the solution is locally stiff
if a perturbation to it remains localized in component number
\emph{and} decays rapidly in time. If every component is locally stiff,
the the solution is globally stiff.

Once strong locality has been established in an ODE system,
individual solution \emph{components}
can be designated as either being (locally) stiff or non-stiff.
For the rest of this paper, when we refer to ``stiff
components'' of the solution, we mean that the components are
locally stiff. A ``non-stiff component'' is one
that evolves on a time scale comparable to the principal time scale.
With this in mind, we can design multirate strategies that handle
stiff and non-stiff components separately. For example,
we expect to be able to integrate all
non-stiff components with large time steps using explicit solvers.
If the number of stiff components is relatively small, we can
use the same explicit method on the stiff components also, but
with much smaller time steps because of stability constraints.
If on the other hand the number of stiff components is fairly
large, we should resort to a fully implicit stiff solver.
%

At this point it is worth comparing our approach with other work in
the literature dealing with problems involving disparate time scales.
Gear and Kevrikidis~\cite{GearKevrikidisSiamJSciComp03} propose their
``projective integration'' method to deal with situations where there
is a \emph{gap} in the spectrum of time scales: the main evolution of a
(stable) solution occurs slowly, with
perturbations decaying much more quickly. For a linear
problem, this corresponds to a situation where the eigenvalues can
be separated into two groups: one set of moderate sized
eigenvalues, and another set with large negative real parts.
Projective integration requires two ODE solvers:
an ``inner'' and an ``outer'' integrator. The idea is
to take many small steps using the inner
solver --- so that the fast modes are damped out, followed by a large
\emph{projective} step with the outer integrator. The process is
then repeated. This method (which is not multirate) is well-suited
to handling problems where many, or all, of the solution components
are rapidly attracted to a slowly varying manifold. Note that
there is no notion of ``locality''
in this approach: the fast modes can potentially be coupled with all
the slow ones. 

Our multirate method also involves an ``inner'' and an
``outer'' integrator. As discussed above, the property
of local stiffness means that a certain subset of the solution
components have much stronger stability than the others.
These stiff components are handled by an ``inner'' integrator
while the non-stiff ones are taken care of using an ``outer''
solver. However, in contrast to the work of Gear and Kevrikidis,
our method is more suited to systems where a small fraction
of the components is stiff at any time during the integration. 
In fact, in terms of the ODE's evolution
in time, one can think of projective integration
as using the inner/outer integrators in ``series'', whereas
our multirate method uses them in ``parallel''.

\subsection[
Strong Locality and Local stiffness in the step flow equations%
]{
Local stiffness for the step flow equations.%
}\label{localSFE}
%
%

We will now show that $\rho_n$, the radius of a bulk step, in (\ref{eqn:sec20:eqnrhon})
is strongly local provided $1 \ll n \ll N$ and hence, away
from the facet and the substrate, any stiffness
that arises is localized. For a bulk step,
the physical origin of the rapid decay comes from the
nature of step interactions in equation (\ref{eqn:sec20:lambda}).
Steps strongly repel each other when they get too close together. 
Consider a configuration
where some of the steps in the bulk are tightly bunched together,
and most of the other steps are widely spaced apart. In this case,
a step strictly (two steps away, at least, from the edge) inside a
bunch is strongly stable, and hence stiff, because small perturbations
in its trajectory are opposed by strong interactions from the
neighboring steps. On the other hand, widely
spaced steps do not experience such large forces, and respond to
perturbations on much slower time scales.
It turns out that these ``step bunching'' 
configurations are quite common in practice and
are produced by the natural time evolution of the system.
In fact, 
the step bunching instability~\cite{IsraeliKandelPRB99,KandelWeeksPRL95,KrugPRB05} 
is a well-studied phenomenon in epitaxial growth, with applications
in quantum dot technology~\cite{KitamuraNishiokaApplPhysLett95} 
and nanolithography~\cite{KeefeUmbachJPhysChemSolids94}.

To analyze the decay of solutions,
the direct approach would be to compute the
Jacobian matrix and analyze its eigenspectrum.
Unfortunately, while the Jacobian
for the system in
(\ref{eqn:sec20:eqnrho1}) -- (\ref{eqn:sec20:eqnrhoN}) can be computed
analytically by linearizing at any fixed set of radii
\Dmath{(\rho_1,\,\rho_2,\, ...,\,\rho_{N-1},\,\rho_N)}, the
expressions involved are very complicated, and do not give much
insight as to why the equations should be stiff. Instead, we present
below a less rigorous calculation, which allows us to relate the
degree of local stiffness to the step spacing. Our approach
is based on the fact that the number of equations, $N$, is
generally rather large, and that the solutions of interest have
a step spacing that is, piecewise, nearly constant.
By this we mean that the step
spacing \Dmath{\;\rho_{n+1}-\rho_n\;} changes slowly with \Dmath{n},
except for a few places where it may change abruptly --- the effect
of these changes is much harder to analyze, and our
method of attack ignores them
since it is only valid far away from these rapid 
transition regions. However, the results of our numerical
calculations indicate that their presence does not invalidate our
analysis.

We begin by considering a configuration of steps which
has a nearly constant step spacing,
and expand the solution in the form $\,\rho_n = \rho_0 +
(n-1)\,\delta\,$ + $\delta^2\,v_n(t) + \dots \,$, where
\Dmath{\rho_0 = O(1)} is a constant, \Dmath{\delta} is
the (constant) leading order step spacing and
$\rho_0 \gg (n-1) \,\delta \gg \delta^2 \,v_n$.
Substituting this
expression into the step flow equations, and ignoring the equations
for the boundary steps (corresponding to \Dmath{n=1\/,\,2\/,N-1\/,} 
and \Dmath{N}),
results in the following leading order equation for the perturbation
$v_n$
\begin{equation} \label{ConstStepSpacingPert}
 \frac{d v_n}{dt}\,  =
 - \frac{3\,\eps}{2\,\delta^4}\,\left(
 v_{n-2} - 4 v_{n-1} + 6 v_n - 4 v_{n+1} + v_{n+2} \right) =
 - \frac{3\,\eps}{2\,\delta^4}\,\left(\Delta^2 \; v\right)_n\/,
\end{equation}
where \Dmath{\Delta} is the discrete Laplacian:
\Dmath{(\Delta \; v)_n = v_{n+1} - 2\,v_n + v_{n-1}}.
To show that equation~(\ref{ConstStepSpacingPert}) has the
property of strong locality, we
consider the solution to the problem with the initial condition
\Dmath{v_n(0) = \nu \delta_{n\/,\,p}}, given by
\begin{equation} \label{ConstStepSpacingIV}
 v_n = \frac{\nu}{2\,\pi} \int_{-\pi}^{\pi}
       \,e^{i\,k\,(n-p) - \sigma(k)\,t}\,dk\/,
 \quad \quad \mbox{where} \quad
 \sigma(k) = \frac{24\,\eps}{\delta^4}\,\sin^4 \frac{k}{2},
\end{equation}
with $\nu \ll 1$. When $|n-p| \gg 1$, $v_n \rightarrow 0$ 
exponentially because the
integrand in (\ref{ConstStepSpacingIV}) is $2\pi$-periodic.
Hence the delta function initial condition remains localized for all $t$.
%

Also, note that equation~(\ref{ConstStepSpacingPert}) has free normal modes given by
\begin{equation} \label{ConstStepSpacingNM}
 v_n = e^{i\,k\,n - \sigma(k)\,t}\/,
\end{equation}
and \Dmath{-\pi < k \leq \pi} is the wave-number. It follows that
the time scales behave like the fourth power of the step
spacing. Hence,
widely spaced steps evolve on a slow time scale whereas
step bunches, which consist of sets of tightly packed steps,
give rise to fast time scales and local stiffness.

\section[
Code Details%
]{
Algorithm Details.%
}\label{sec04:code}
\subsection[
Algorithm Details%
]{
Algorithm Overview%
}\label{algorithm}
%

The goal of our method is to efficiently solve a system of locally
coupled ODEs where only a few of the components are stiff. A standard
explicit integrator would take small time steps for all components of
the solution. In contrast, our multirate method takes 
large steps for the non-stiff components,
and small steps for the stiff ones.


The algorithm starts by taking an explicit, global time step, say:
from $t_n\/$ to $t_{n+1}\/$. An embedded formula is then 
used to obtain an estimate of the Local
Truncation Error (LTE) for each component of the solution. In general,
some of the LTEs will be unacceptably large (because 
the associated solution components are stiff), 
while others will have acceptable sizes. The algorithm
checks if the components with acceptable LTEs satisfy
the preset tolerance levels. If they do not, the step size
is reduced and another global time step is attempted. If they do,
a second round of integration is performed to correct the components
with large LTEs. Hence, the algorithm is as follows:

\begin{PWFalgorithm}~~\newline
\begin{itemize}
\item[1.] Take a step from $t_n$ to $t_{n+1}$.
\item[2.] \vspace*{-1mm} 
  Let \Dmath{e_i} be the LTE for the \Dmath{i}-th component, let
  \Dmath{\mu} be the $P^{th}$ percentile of all the LTEs, and let \Dmath{tol_i} be
  the required error tolerance for $i$-th component. For example, if $P=10$, then
  $90 \%$ of the errors are larger than $\mu$. If $P=50$, $\mu$ is the median.

  \item[3.] 
  For some real number \Dmath{k}, flag all the components whose
  LTEs are greater than \Dmath{10^k \mu} (our code uses \Dmath{k=2}) as
  being possibly stiff.

  \item[4.] 
  Check if the unflagged solution components satisfy the tolerance
  requirements, i.e. $M \equiv Max(|e_i/tol_i|) < 1$ where $Max()$
  is taken over all unflagged components only.
  \item[5.~] \begin{itemize}
   \item[(a)] If they do, the step is successful. 
  The step size is increased using the formula 
  (\ref{stepsizecontrol}) with $M^{-1} > 1$ as the ratio of errors.
  \item[(b)] Perform a second integration 
  to correct the flagged components (see Figure \ref{fig4}). 
  \item[(c)] Increase the step size according
  to formula (\ref{stepsizecontrol}). Increase $n$ and go back to step 1.
  \end{itemize}
  \item[6.~] \begin{itemize}
    \item[(a)] If they do not, the step is not successful. Reduce the step size 
    according to (\ref{stepsizecontrol}) with $M^{-1} < 1$ as the ratio of errors.
    \item[(b)] Do not increase $n$ and go back to step 1.
  \end{itemize}
\end{itemize}
\label{algorithm2}
\end{PWFalgorithm}

The second integration is basically done only for the
stiff components, and it involves many small sub-steps within
the interval \Dmath{[t_n,t_{n+1}]} to ensure stability.
Although this second round of integration takes a
large number of steps, it only needs to be done for a
small subset of the solution components.

\RRRfigOLD{\EPfg{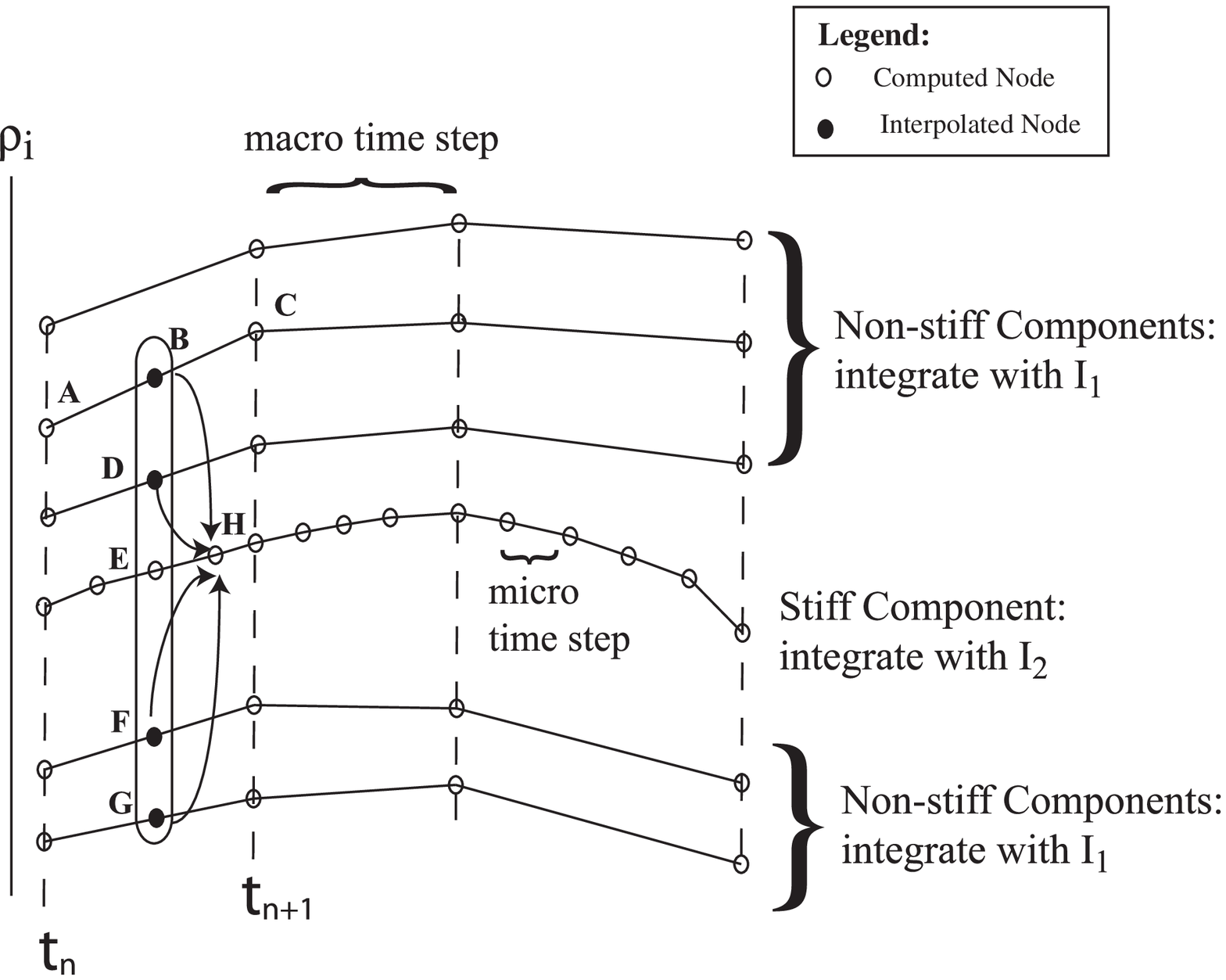,width=5in}}{
 \caption{
 Schematic showing the two phases of integration with the time-steppers
 \Dmath{I_1} and  \Dmath{I_2}. Large time steps are taken
 using $I_1$ in the slowly varying bulk while smaller time
 steps are taken for stiff components using $I_2$.
 Here, point B is
 interpolated from points A and C, with D, F, and G obtained in a
 similar fashion. Because equations
 (\ref{eqn:sec20:eqnrho1} -- \ref{eqn:sec20:eqnrhoN})
 are pentadiagonal, all the points B, D, E, F, and G are needed to
 to compute point H. Finally, we point out that in practice the
 stiff components are not isolated, but appear in bunches.}
 \label{fig4}}

To perform the second integration, values for the non-stiff components
at all times $t \in (t_n,t_{n+1})$ along the boundaries of any stiff 
set of steps are needed --- see Figure
\ref{fig4}. For example, in the case of a pentadiagonal system, the
values of two non-stiff components are needed on each side of a stiff
region. One way to generate dense output from the
non-stiff components, between \Dmath{t_n} and \Dmath{t_{n+1}}, is
through interpolation. In this paper we use cubic interpolatants
which are generated by using the intermediate stage function evaluations
in a Runge-Kutta method. In Section \ref{interpolation},
we give more details on the construction of these interpolants and demonstrate that
cubic interpolation is consistent with a multirate method
that is globally fourth order.

Once solution components have been flagged as requiring re-integration, the local coupling
means that some of the non-stiff components may also have to be re-integrated. 
Because (\ref{eqn:sec20:eqnrho1}) -- (\ref{eqn:sec20:eqnrhoN})
is a pentadiagonal system of equations, if only $r_m(t)$ 
and $r_{m+2}(t)$ are stiff
components with large LTEs, then all three of the components $r_m(t),r_{m+1}(t)$ and $r_{m+2}(t)$ must
be re-integrated as a set, using the dense output from
$r_{m-2}(t)$, $r_{m-1}(t)$, $r_{m+3}(t)$ and $r_{m+4}(t)$ as `boundary
conditions'. Hence, the algorithm is slightly wasteful in that although $r_{m+1}(t)$ was deemed
accurate enough, it still had to be integrated for a second time.

Note that our algorithm uses the LTE in 
a different way from conventional
embedded RK methods: Instead of immediately
scaling the time step if the smallest LTE is greater 
than the tolerance level, 
we make a note of \textit{which} components 
had the largest LTEs by analyzing their distribution:
it might not be efficient
to retake the time step for every component, if 
only a few of them are inaccurate. 
The largest LTEs (in the sense of being larger
than $10^k \mu$ in Algorithm \ref{algorithm2}) are discarded, and
then the time step scaled according to the 
largest of the remaining errors. 
Hence, we get a larger time step for a majority
of the solution components, and the way that 
this time step is adjusted throughout the course
of the integration is not affected by the presence 
of a few either rapidly varying or stiff components.

For the rest of this paper, we will call the first 
time stepper $I_1$ (used to generate the LTEs
in the first place), and the second time stepper $I_2$ 
(used to re-integate stiff components with
large LTEs). In general, $I_1$ and $I_2$ do not have 
to be the same method, or of the same
order, but $I_1$ has to be
able to generate estimates of the Local Truncation Error. 
In our code, $I_1$ is a 
Cash-Karp Runge-Kutta Formula \cite{Press} and $I_2$ 
depends on the solution component: if the 
re-integration involves the innermost
step (see Section \ref{treatment}), we take $I_2$ 
to be a Simple-Euler routine
which adjusts its step size by step doubling; otherwise $I_2 = I_1$.
In other words, there are two possibilities which 
can arise when performing the
re-integration with $I_2$:
\begin{enumerate}
\item The re-integration involves solution components 
which include the innermost step. In the following analysis,
we will assume this is $\rho_1$.
\label{case1}
\item The re-integration does not involve the innermost step. 
\label{case2}
\end{enumerate}
The reason to distinguish between these two 
cases is that (\ref{case1}) will involve integration
of singular trajectories (see equation (\ref{square_root})), 
but in general, (\ref{case2}) will not. 
Using step doubling in (\ref{case1}) is a 
fairly crude way of adjusting the time step. 
However, resorting to an embedded formula
is not possible when
(\ref{stepsizecontrol}) breaks down, so
using step doubling to monitor the quality of the 
solution is reasonable in this case.

\subsection{Treatment of Singular Collapse of Top Step} 
\label{treatment}
From equation (\ref{square_root}), we have seen that $\rho_1 \rightarrow 0$
in a singular fashion, causing problems for standard
high order integrators. Our treatment
uses a Simple Euler method whenever the re-integration
of $\rho_1$ and neighboring singular components is involved: this is ``optimal'' in the sense
that Simple Euler produces results that have the same accuracy as higher
order methods (due to the singular nature of $\rho_1$) but
is computationally cheaper.
Furthermore, we are able
to extract the time of collapse, $t_1$, using linear
interpolation, which is consistent with Simple Euler's order
of accuracy.

Our method involves solving for 
$\rho_1^2,\rho_2,...\rho_N$ instead of $\rho_1,\rho_2,...,\rho_N$. Note that
from (\ref{square_root}),
\begin{equation}
\rho_1 ^2(t) \sim C_3 (\tau_1 - t) + C_4 (\tau_1 - t)^{3/2} \ln(\tau_1 - t) \label{rho_sq}
\end{equation}

as $t \rightarrow \tau_1$ (the collapse time) for some constants $C_3$ and $C_4$,
which means that $\rho_1^2$ has exactly one derivative at $\tau_1$.
Our main reason for solving for $\rho_1^2$, instead
of $\rho_1$, is not to improve accuracy, but rather to 
enable the algorithm to `step through' 
the singularity at $t=\tau_1$, and use linear
interpolation to obtain $\tau_1$, the time of collapse 
of the innermost step. 

Taking square roots to recover $\rho_1$ will
will result in a drastic loss in accuracy near 
$\tau_1$. At time $t$ close to $\tau_1$, 
consider taking a time step of size
$\Delta t$ with component $\rho_1^2$ using Simple 
Euler. Let $\rho^2_{\textrm{exact}}(t+\Delta t)$ be the result of taking this
time step using a `perfect' integrator, producing 
the exact solution at $t+\Delta t$, given $\rho_1^2(t)$. Then, since
the truncation error in Simple Euler is $O\left(\Delta t^2 
\dvn{\rho_1^2}{t}{2}\right) = O(\Delta t^{3/2} \ln \Delta t)$ from (\ref{rho_sq}), we have

\begin{eqnarray}
|\rho_1^2(t+\Delta t) - \rho^2_{\textrm{exact}}(t+\Delta t)| &=& O(\Delta t^{3/2} \ln \Delta t), \\
\Rightarrow \rho_1 &=& \rho_{\textrm{exact}} \left( 1+ \frac{O(\Delta t^{3/2} 
\ln \Delta t))}{\rho_{\textrm{exact}}^2}\right)^{1/2}.
\end{eqnarray}
Therefore, if $\rho^2_{\textrm{exact}} \gg O(\Delta t^{3/2} 
\ln \Delta t)$ ($t$ sufficiently far away from the singularity)
then the LTE for $\rho_1$, $|\rho_1 - \rho_{\textrm{exact}}|$, 
is $O(\Delta t^{3/2} \ln \Delta t)$. However, if
$\rho^2_{\textrm{exact}} \ll O(\Delta t^{3/2} \ln \Delta t)$
($t$ is very close to $t_1$), then the LTE for 
$\rho_1$ is $O(\Delta t^{3/4} (\ln \Delta t)^{1/2})$, 
which is not a big improvement over (\ref{eq100}).
Note that these estimates for the LTE are 
independent of the order of $I_2$.
When $I_2$ has `overstepped' $\tau_1$
resulting in $\rho_1^2(t_m) > 0$ and $\rho_1^2(t_{m+1}) < 0$ for 
times $t = t_m, t_{m+1} \equiv t_m + \Delta t$ ($t_m < \tau_1 < t_{m+1}$),
we set
\begin{equation}
\tau_1 \approx \frac{ \rho_1^2(t_{m+1}) t_{m+1} }{ \rho_1^2(t_m)-\rho_1^2(t_{m+1}) }
- \frac{\rho_1^2(t_{m+1}) t_m}{\rho_1^2(t_m)-\rho_1^2(t_{m+1})},
\end{equation}
as an approximation to the collapse time. Once $\rho_1$ has collapsed 
at $\tau_1$, it is removed from the system
(\ref{eqn:sec20:eqnrho1}) -- (\ref{eqn:sec20:eqnrhoN}),
the number of equations drops by one, and $\rho_2^2(t)$
replaces $\rho_1^2(t)$ as the new top step.

\subsection{Interpolation}
\label{interpolation}

The key to making our multirate method high order lies in
the ability to generate dense output 
from the non-stiff components with high accuracy.
One way to generate dense output
is to use interpolation\footnote{
Although we use the word ``interpolation'' to describe
a method to construct a continuous function
between the two time
points $t_n$ and $t_{n+1}$, the function that we derive does
not actually pass through $t_{n+1}$. Hence 
strictly speaking, it is not an interpolant, though we will
continue to refer to these approximating functions as
``interpolants'' for convenience.}. For some integrators,
such as Backward Differentiation Formulae (BDF), 
is it obvious how to derive an interpolant that is
consistent in order with the underlying integrator -- BDF
use extrapolation to advance the solution in time.
For other integrators, such as Runge-Kutta schemes, constructing
the interpolant is less obvious and this is the focus of the
section. Note that
we need to generate interpolants
\emph{during run time} using only the function evaluations
that have already been computed by the integrator 
within each time step.
The extra constraint of generating the interpolants during run time
adds a non-trivial complication
to the ``traditional'' interpolation
problem, which has been studied extensively 
\cite{EnrightJacksonACMTransactions86,HornSIAM83}.
Having successfully integrated the non-stiff 
components from $t_n$ to $t_{n+1}$, we have
the point values $y_n, y_{n+1}$ and the derivative 
$y_n'$ at our disposal to construct
the interpolant between $t_n$ and $t_{n+1}$. We do \emph{not}
have information about the derivative $y_{n+1}'$.
However, because we are constructing 
these interpolants during run time, we are at liberty to
use the intermediate function evalulations 
inherent in the application of Runge-Kutta: this
is valuable information that is not 
usually available in the traditional interpolation problem.

We have seen that our method performs integration in two phases:
first we integrate a large number of non-stiff components, then we integrate
a small number of stiff ones. Ideally, we would
like the two integrations to have the same order.
This is only possible if the interpolation of the non-stiff components
is of a sufficiently high order -- otherwise large interpolation errors will 
contaminate the accuracy in the stiff components. In the following paragraphs,
when we generate dense output between points $t$ and $t+\Delta t$,
we define an interpolant to be order $m$ when the interpolation error
is $O(\Delta t^{m+1})$.

Let us assume that our inner and outer solvers are both $n^{th}$ order.
First of all, let us calculate $m$ in terms of $n$
if we want our method to be globally $n^{th}$ order. 
For simplicity, we assume in this calculation that the effects of round-off error are negligible.
Consider the ODE system 
\begin{equation}
\mathbf{y}' = \mathbf{F}(\mathbf{y},t).
\end{equation}
Let us assume that we have taken a macro step of size 
$\Delta t$ and advanced the non-stiff
components $\mathbf{y}_r(t)$ successfully from time $t$ 
to $t+\Delta t$. Also assume
that we have taken $N$ micro steps
of size $\Delta t_i$, $i=1,2,...,N$, for the stiff 
components $\mathbf{y}_s(t)$ so that $\sum_{i=1}^N \Delta t_i = \Delta t$.
After taking these $N$ steps, the total error in 
$\mathbf{y}_s(t)$ will be
\begin{equation}
O\left( \sum_{i=1}^N \Delta t_i^{n+1}\right) + O\left(\sum_{i=1}^N 
\Delta t_i \Delta t^{m+1} \right).
\label{toterr}
\end{equation}
The first term is the sum of the Local Truncation 
Errors caused by taking $N$
steps each of size $\Delta t_i$. The second term is 
the sum of the interpolation errors:
note that to advance $\mathbf{y}_s$, an evaluation of 
$\mathbf{F}$ in between $t$
and $t+\Delta t$, in general, is required 
for the non-stiff neighbours of $\mathbf{y}_s$ 
and this incurs an interpolation error of size 
$O(\Delta t^{m+1})$. Therefore
the error in $\mathbf{y}_s'$ is also $O(\Delta t^{m+1})$ and the 
error in $\mathbf{y}_s$ is $O(\Delta t_i \Delta t^{m+1})$.
It is clear that (\ref{toterr}) simplifies to
\begin{equation}
O\left( \Delta t^{n+1}, \Delta t^{m+2} \right),
\end{equation}
and so for our multirate method to be globally $n^{th}$ order, we require $m=n-1$ -- that
is, we can afford for the order of the interpolation to be one less than the order of the integrator.
For example, if our integrator is fourth order ($n=4$), we need to be able to construct
cubic interpolants during run time. If our integrator is second order, then linear interpolation
should be sufficient -- as observed in \cite{SavcencoHundsdorferBIT06}. 

We will now illustrate how these interpolants are constructed by taking the classical 4th order
(non-adaptive) Runge-Kutta formula as an example and applying it to the autonomous ODE
system $y'=f(y)$:

\begin{table}[htbp]
\caption{Coefficients in classical RK4}
\label{rktab}
\begin{center}
\begin{tabular}{l|llll}
$a_1 = 0$ & $b_{11} = 0$ &&& \\
$a_2 = \frac{1}{2}$ & $b_{21} = \frac{1}{2}$ & $b_{22} = 0$ && \\
$a_3 = \frac{1}{2}$ & $b_{31} = 0$ & $b_{32} = \frac{1}{2}$ & $b_{33} = 0$ & \\
$a_4 = 1$ & $b_{41} = 0$ & $b_{42} = 0$ & $b_{43} = 1$ & $b_{44} = 0$\\
\hline
& $c_1 = \frac{1}{6}$ & $c_2 = \frac{1}{3}$ & $c_3 = \frac{1}{3}$ & $c_4 = \frac{1}{6}$ \\  
\end{tabular}
\end{center}
\end{table}
The solution $y_n$ is advanced to $y_{n+1}$ through
\begin{equation}
\label{rkstep}
y_{n+1} = y_n + \sum_{i=1}^4 c_i k_i,
\end{equation}
where
\begin{equation}
k_i = \Delta t f(y_n + \sum_{j=1}^{i-1} b_{ij} k_j ),
\end{equation}
for $i=1,2,3,4$. It is easy to show that equation (\ref{rkstep}) implies that 
\begin{equation}
\label{ynp1}
y_{n+1} = y_n + \Delta t y_n' + \frac{\Delta t^2}{2}y_n'' + \frac{\Delta t^3}{3!}y_n''' + \frac{\Delta t^4}{4!}y_n^{(4)} + O(\Delta t^5).
\end{equation}
Let us try to construct a quartic interpolant. Equation (\ref{ynp1}) motivates us
to write this in the form
\begin{equation}
y(\chi \Delta t) = y_n + (\chi \Delta t) y_n' + \frac{(\chi \Delta t)^2}{2}y_n'' + \frac{(\chi \Delta t)^3}{3!}y_n'''
+ \frac{(\chi \Delta t)^4}{4!} y_n^{(4)} + O(\Delta t^5),
\end{equation}
where $0 \leq \chi \leq 1$. The problem is now to evaluate the derivatives $y_n^{(m)}$ in terms
of the intermediate stage function evaluations $k_i$. This is done by (i) noting that
\begin{equation} \label{yn}
\begin{array}{lll}
y_n'  & = & f_n \\
y_n'' & = & f_n' f_n \\
y_n''' & = & f_n {f_n'}^2 + f_n^2 f_n'' \\ 
y_n^{(4)} & = & f_n {f_n'}^3 + 4 f_n^2 f_n' f_n'' + f_n^3 f_n'''
\end{array}
\end{equation}
where $f_n \equiv f(y_n)$ and similarly with $f_n', f_n'',...$, 
and (ii) expanding $k_i$ in Taylor series:
\begin{equation} \label{ki}
\begin{array}{lll}
k_1 &=& \Delta tf_n, \\
k_2 &=& \Delta tf_n + \frac{\Delta t^2}{2} f_n'f_n + 
\frac{\Delta t^3}{8} f_n^2 f_n'' + 
\frac{\Delta t^4}{48} f_n^3 f_n''' + O(\Delta t^5), \\
k_3 &=& \Delta tf_n + \frac{\Delta t^2}{2}f_n'f_n + 
\frac{\Delta t^3}{8} (2 f_n {f_n'}^2 + f_n^2 f_n'') \\
&& + \frac{\Delta t^4}{48} ( 9 f_n^2 f_n' f_n'' +
f_n^3 f_n''') + O(\Delta t^5), \\
k_4 &=& \Delta tf_n + \Delta t^2 f_n'f_n + 
\frac{\Delta t^3}{2} (f_n {f_n'}^2 + f_n^2 f_n'') \\
&& + \frac{\Delta t^4}{24} ( 6f_nf_n'^3 + 15 f_n^2 f_n'f_n'' + 
4f_n^3 f_n''') + O(\Delta t^5).
\end{array}
\end{equation}
A natural way to compute the $y_n^{(m)}$ would be to find the 7 terms $f_n$, $f_n' f_n$, $f_n {f_n'}^2$, $f_n^2 f_n''$,
$f_n {f_n'}^3$,  $f_n^2 f_n' f_n''$ and $f_n^3 f_n'''$ in terms of the $k_i$ from (\ref{ki}) and then use them in (\ref{yn}).
However, this is not possible because (\ref{ki}) becomes a system of 4 linear equations in 7 unknowns. We must therefore be
a little less ambitious. In light of our previous comments on interpolation, we seek a \emph{cubic} interpolant in the form
\begin{equation}
y(\chi \Delta t) = y_n + (\chi \Delta t) y_n' + \frac{(\chi \Delta t)^2}{2}y_n'' + \frac{(\chi \Delta t)^3}{3!}y_n'''
+ O(\Delta t^4).
\end{equation}
Since we do not need the $y_n^{(4)}$ term, constructing $y_n', y_n''$ and $y_n'''$ now requires only 
$f_n$, $f_n' f_n$, $f_n {f_n'}^2$ and $f_n^2 f_n''$. Ignoring the $O(\Delta t^4)$ terms, (\ref{ki}) constitutes
4 equations in 4 unknowns. Solving in terms of the $k_i$ yields
\begin{eqnarray}
f_n &=& \frac{k_1}{\Delta t}, \\
f_n f_n' &=& \frac{1}{\Delta t^2} \left( -3k_1 + 2k_2 + 2k_3 - k_4  \right),\\
f_n^2 f_n'' &=& \frac{1}{\Delta t^3} \left( 4 k_1 -8 k_3 + 4k_4 \right), \\
f_n {f_n'}^2 &=& \frac{1}{\Delta t^3} \left( -4k_2 + 4k_3 \right),
\end{eqnarray}
so that the cubic interpolant is
\begin{equation}
y(\chi \Delta t) = y_n + \chi k_1 + \frac{\chi^2}{2}\left( -3k_1 + 2k_2+2k_3-k_4\right) + \frac{2\chi^3}{3}\left(k_1-k_2-k_3+k_4\right).
\end{equation}
We now turn our attention to Embedded Runge-Kutta Methods. Let us focus on the Cash-Karp formula
in \cite{Press} which has the tableau

\begin{table}[htbp]
\begin{center}
\caption{Coefficients in the Cash-Karp 4-5 formula}
\label{cktab}
\begin{tabular}{l|llllll}
$0$     & 0                      &&&&& \\
$\frac{1}{5}$     & $\frac{1}{5}$ & 0 &&&& \\
$\frac{3}{10}$    & $\frac{3}{40}$ & $\frac{9}{40}$ & 0 &&& \\
$\frac{3}{5}$     & $\frac{3}{10}$ & $-\frac{9}{10}$ & $\frac{6}{5}$ & 0 && \\
$1$               & $-\frac{11}{54}$ & $\frac{5}{2}$ & $-\frac{70}{27}$ & $\frac{35}{27}$ & 0 & \\
$\frac{7}{8}$     & $\frac{1631}{55296}$ & $\frac{175}{512}$ & $\frac{575}{13824}$&
$\frac{44275}{110592}$ & $\frac{253}{4096}$ & 0 \\ 
\hline
& $\frac{2825}{27648}$ & $0$ & $\frac{18575}{48384}$ & $\frac{13525}{55296}$ & $\frac{277}{14336}$ & $\frac{1}{4}$\\
& $\frac{37}{378}$ & $0$ & $\frac{250}{621}$ & $\frac{125}{594}$ & $0$ & $\frac{512}{1771}$
\end{tabular}
\end{center}
\end{table}

The analogue to (\ref{ki}) is
\begin{equation}\label{ki2}
\begin{array}{lll}
k_1 &=& \Delta t f_n, \\
k_2 &=& \Delta t f_n + \frac{\Delta t^2}{5}f_n' f_n + 
\frac{\Delta t^3}{50} f_n^2 f_n'' + 
\frac{\Delta t^4}{750} f_n^3 f_n''' + O(\Delta t^5), \\
k_3 &=& \Delta t f_n + \frac{3\Delta t^2}{10}f_n' f_n + 
\frac{9\Delta t^3}{200} (f_n {f_n'}^2 + f_n^2 f_n'') \\
&& + \frac{9\Delta t^4}{2000} ( 4f_n^2 f_n' f_n'' +
f_n^3 f_n''') + O(\Delta t^5), \\
k_4 &=& \Delta t f_n + \frac{3\Delta t^2}{5}f_n' f_n + 
\frac{9\Delta t^3}{50} (f_n {f_n'}^2 + f_n^2 f_n'') \\
&& + \frac{9\Delta t^4}{500} ( 3f_nf_n'^3 + 8f_n^2 f_n'f_n'' +
2f_n^3 f_n''') + O(\Delta t^5), \\
k_5 &=& \Delta tf_n + \Delta t^2 f_n' f_n + 
\frac{\Delta t^3}{2} (f_n {f_n'}^2 + f_n^2 f_n'') \\ 
&& + \frac{\Delta t^4}{60} ( 7f_nf_n'^3 + 40 f_n^2 f_n' f_n'' + 10f_n^3 f_n''' )
+ O(\Delta t^5), \\
k_6 &=& \Delta tf_n + \frac{7\Delta t^2}{8}f_n' f_n + 
\frac{49\Delta t^3}{128}(f_n {f_n'}^2 + f_n^2 f_n'') \\
&& + \frac{7\Delta t^4}{3072} ( 46f_nf_n'^3 + 196 f_n^2 f_n'f_n''
+ 49 f_n^3 f_n''') + O(\Delta t^5).
\end{array}
\end{equation}
Cash-Karp 4-5 is formally 4th order, so again, it is sufficient to interpolate with
cubic polynomials. However, one would expect that since we have made two extra evaluations of $f(y)$, it would
be possible to construct interpolants which have higher order. At first, this possibility
seems promising. We note that in (\ref{ki2}), with the exception of the $k_2$ equation, 
the terms $f_n^2 f_n' f_n''$ and $f_n^3 f_n'''$ always appear together as $(4 f_n^2 f_n' f_n'' + f_n^3 f_n''')$ and
$f_n {f_n'}^2$, $f_n^2 f_n''$ always appear together as $(f_n {f_n'}^2 + f_n^2 f_n'')$.
The equations for $k_1,k_3,k_4,k_5$ and $k_6$ therefore give 5 equations in 5 unknowns, $f_n$, $f_n' f_n$,
$(f_n {f_n'}^2 + f_n^2 f_n'')$, $f_n f_n'^3$ and $(4 f_n^2 f_n' f_n'' + f_n^3 f_n''')$. Unfortunately, the
matrix of the resulting linear system has a zero determinant and the equation
\begin{equation}
\label{mateqn}
\begin{pmatrix}
\begin{array}{lllll}
\Delta t & 0 & 0 & 0 & 0 \\
\Delta t & \frac{3\Delta t^2}{10} & \frac{9\Delta t^3}{200} & 0 & \frac{9\Delta t^4}{2000} \\
\Delta t & \frac{3\Delta t^2}{5}  & \frac{9\Delta t^3}{50}  & \frac{27\Delta t^4}{500} & \frac{9\Delta t^4}{250} \\
\Delta t & \Delta t^2             & \frac{\Delta t^3}{2}    & \frac{7\Delta t^4}{60}    & \frac{\Delta t^4}{6}   \\
\Delta t & \frac{7\Delta t^2}{8}  & \frac{49\Delta t^3}{128}& \frac{161\Delta t^4}{1536}& \frac{343\Delta t^4}{3072}
\end{array}
\end{pmatrix}
\begin{pmatrix}
\begin{array}{c}
f_n \\
f_n' f_n \\
f_n {f_n'}^2 + f_n^2 f_n'' \\
f_n {f_n'}^3 \\
4 f_n^2 f_n' f_n'' + f_n^3 f_n'''
\end{array}
\end{pmatrix} = 
\begin{pmatrix}
\begin{array}{c}
k_1 \\
k_3 \\
k_4 \\
k_5 \\
k_6
\end{array}
\end{pmatrix} 
\end{equation}
does not have a unique solution.

Going through the same process with Fehlberg's
pair does not improve the situation. It is not clear to us at this point whether the inability to 
build quartic interpolatants using the intermediate stage function evaluations
is symptomatic of all RK45 pairs, or if it is possible to find Runge-Kutta families for which
the equivalent of (\ref{mateqn}) is uniquely solvable. Although having quartic interpolants
is not necessary for our multirate method to be globally fourth order, these interpolants
-- if they can be constructed -- could be used to check the accuracy of the cubic interpolants
in the same way that the fifth order RK formula is used to check the 
values predicted by the fourth
order one. If the interpolation is deemed too inaccurate, the integration of the
non-stiff components would have to be performed again using a smaller step size.

Construction of the cubic interpolant in Cash-Karp 45 is now fairly straightforward and
we follow the same procedure as for classical RK4. There is now more than one possible
cubic polynomial, depending on which $k_i$ to use in (\ref{ki2}).
Using the equations involving $k_1$, $k_4$ and $k_5$, we have
\begin{equation}
y(\chi \Delta t) = y_n + \chi k_1 + \frac{\chi^2}{2}\left( -\frac{8}{3}k_1 + \frac{25}{6} k_4
-\frac{3}{2}k_5 \right) + \frac{\chi^3}{6}\left(\frac{10}{3}k_1 - \frac{25}{3}k_4
+ 5 k_5 \right).
\end{equation}
This is the interpolant used in our multirate code to generate the results in Section \ref{sec05:implementation}.

\section{Validation and Results}
\label{sec05:implementation}
Here we validate our code with different tests, each of which examines a
particular aspect of the integration.

\subsection{Validation}
\subsubsection{Collapse Times}
To test the code's ability to handle singular collapses, we used it to solve
the (uncoupled) ODE system
\begin{equation}
\dot{r}_i = -1/r_i \label{eq1}
\end{equation}
for $i=1,2,...,N$ with initial condition $r_i(0) = i$.

The solution to this set of ODEs is $r_i(t) = \sqrt{i^2 - 2t}$. Note that 
the solution has the same leading order 
singular behaviour at the collapse times
$t_i = i^2/2$ as equations (\ref{eqn:sec20:eqnrho1})-(\ref{eqn:sec20:eqnrhoN}). 
A second, more challenging, model problem is
\begin{equation}
\dot{r}_i = -1/r_i^2, \label{eq2}
\end{equation}
with the same initial condition as before. The
collapse times in this system take the form $t_i = i^3/3$ because the exact solutions are
$r_i(t) = (i^3-3t)^{1/3}$: the solution near the collapse times in this case are even
steeper than in (\ref{eq1}). The results in Table \ref{table1} show that
our code is able to capture the collapse times to 4-5 significant digits of accuracy.

\vspace{0.4cm}
\begin{table}[htbp]
\caption{First 5 collapse times for model systems (\ref{eq1}) and (\ref{eq2}). Both exact
and numerical values are shown. The relative and absolute tolerances were $10^{-6}$ and $10^{-8}$
respectively.}
\vspace{0.4cm}
\begin{center}
\begin{tabular}{c p{1.2in} p{1in} p{1.2in} p{0.6in}}
\hline
Collapse time & \textit{eqn} (\ref{eq1}) & \textit{eqn} (\ref{eq1})
& \textit{eqn} (\ref{eq2}) & \textit{eqn} (\ref{eq2}) \\
& \emph{numerical} & \emph{exact} & \emph{numerical} & \emph{exact} \\
\hline
$t_1$ & \textbf{0.500000000} & 0.50 &  \textbf{0.3333}4710950
   &  0.333\dots \\
$t_2$ & \textbf{2.0000}2007041 & 2.00 &  \textbf{2.666}73866493
   &  2.666\dots \\
$t_3$ & \textbf{4.5000}6868256 & 4.50 &  \textbf{9.000}20498167
   &  9.000\dots \\
$t_4$ & \textbf{8.000}13388581 & 8.00 & \textbf{21.333}76504303
   & 21.333\dots \\
$t_5$ & \textbf{12.500}2133304 & 12.50 & \textbf{41.66}743834183
   & 41.666\dots \\
\hline
\label{table1}
\end{tabular}
\end{center}
\end{table}

To test the accuracy of the collapse times generated from the full set of
step flow equations (\ref{eqn:sec20:eqnrho1}) - (\ref{eqn:sec20:eqnrhoN}), 
we used fixed step Simple Euler with $\Delta t = 10^{-6}$
and linear interpolation to obtain a set of reference
collapse times. Since this method of integration is computationally
expensive, we initialized the simulation with only $N=15$ layers. 
These times were compared with data generated from the full multirate
code with high order adaptive time stepping.
The results are shown in Table \ref{table2}.

\vspace{0.4cm}
\begin{table}[htbp]
\caption{First five collapse times shown for a Diffusion Limited system
(obtained by setting $m_1 = 1$, $m_2 = 0$ in equation (\ref{eqn:sec20:Deltan})
and $\gamma = 1$ in (\ref{eqn:sec20:eqnrho1})-(\ref{eqn:sec20:eqnrhoN})), 
with $\varepsilon=0.01$.
The initial condition was a 15 layer profile with unit spacing.}\label{table2}
\vspace{0.4cm}
\begin{center}
\begin{tabular}{c p{1.5in} p{1.5in}}
\hline
Collapse time & Reference Solution & Multirate Solution\\
\hline
$t_1$ & 0.540289230794 & \textbf{0.540}305641980 \\
$t_2$ & 5.100219762927 & \textbf{5.1002}84674837 \\
$t_3$ & 21.036583847637 & \textbf{21.036}757035271 \\
$t_4$ & 59.481455149416 & \textbf{59.481}830949331 \\
$t_5$ & 135.366866973862 & \textbf{135.36}7562952919 \\
\hline
\end{tabular}
\end{center}
\end{table}

\subsubsection{Convergence and error analysis}
\label{subsec:error-spread}
\begin{figure}
\begin{center}
\epsfig{file=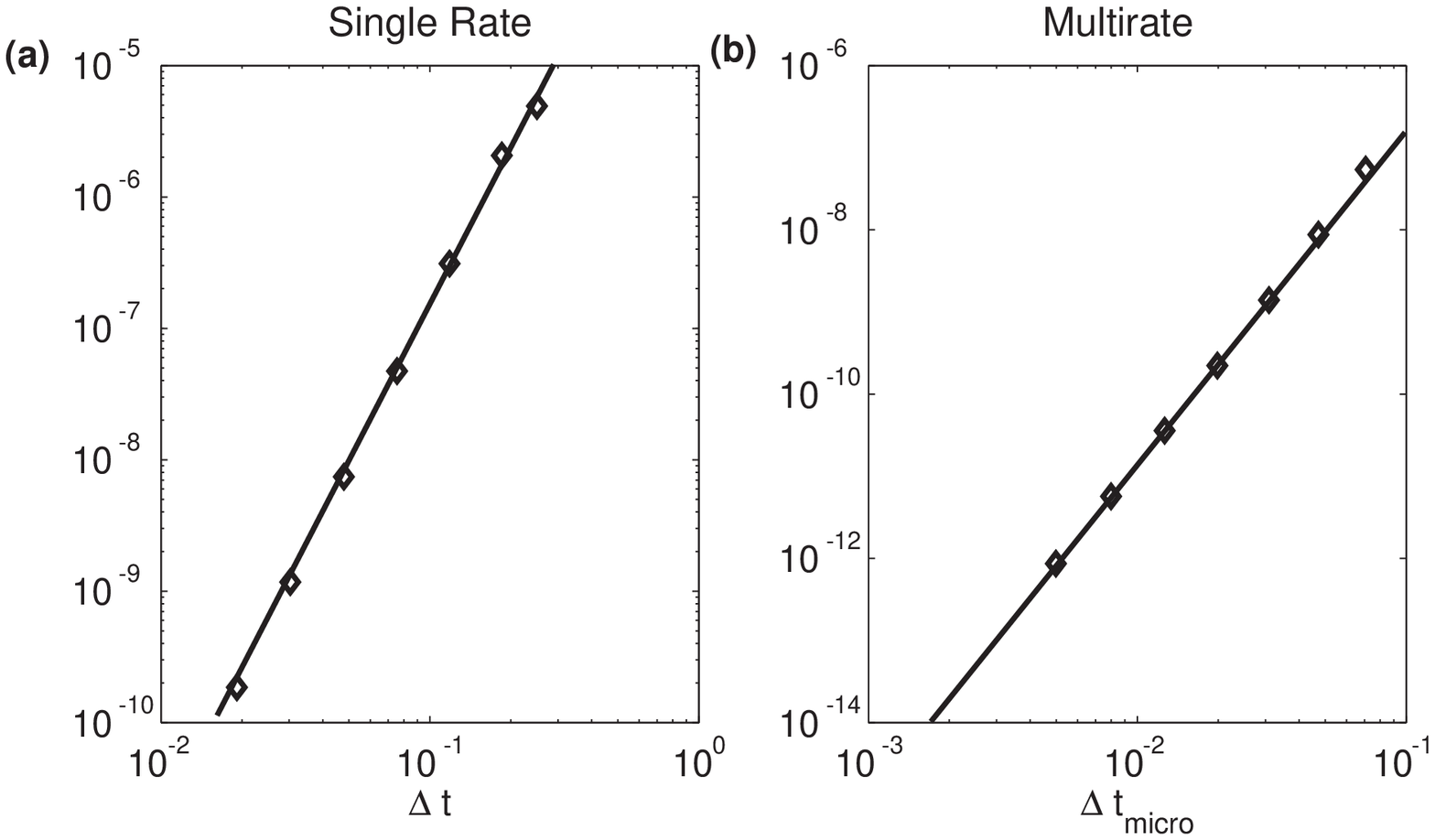,scale=0.7}
\end{center}
\caption{Fourth order convergence was obtained from
solving a wave equation discretized using the method
of lines. The domain was $-25 \leq x \leq 25$ and $N=401$
equations were solved. The final
integration time was $T=20$. Different (average)
time step sizes were obtained by changing the tolerance level of the code.
(a) Fourth order convergence of single rate Cash-Karp Runge Kutta. (b)
Fourth order convergence of multirate Cash-Karp Runge-Kutta.
The parameters used in Algorithm \ref{algorithm2} 
were $k=-6$, $P=30$. The number of macrosteps
used by the algorithm in every integration was 22.}
\label{fig:convergence}
\end{figure}

We used our multirate code to solve the simple wave equation $u_t + u_x = 0$.
By discretizing using the Method of Lines and one-sided (``upwind'') differences in space,
we obtained $N$ coupled Ordinary Differential Equations which were solved
using the initial condition $u(x,0) = e^{-(x+10)^2}$ and periodic
boundary conditions.\footnote{These conditions were implemented by ensuring
the domain of solution was large enough so that the use of periodic boundary
conditions did not introduce any significant errors.} 
Unlike the axisymmetric step equations, this discretization of the
wave equation does not have any time singularities and we expect
fourth order convergence for every component.
This is confirmed in Figure \ref{fig:convergence}.
Although $\Delta t$ and $\Delta t_{\textrm{micro}}$ constantly change because our
algorithm uses adaptive step size control, we take $\Delta t \sim T/N$ 
as a measure of the average
step size where $T$ is the final integration time. For the
single rate method, $N$ is the number
of steps and for the multirate method, $N$ is the total number of micro-steps.
In the multirate code, a maximum macro-stepsize $\Delta t_{\textrm{macro}} = 1$
was imposed, and all but the first and final macro-steps had size 1.
The exact solution comes 
from solving the linear system of ODEs exactly
using an eigenfunction expansion.
As the integration progresses, the fraction of components 
that has to be reintegrated increases gradually
as the solution broadens and its amplitude decreases: see Figure \ref{fig:wave}.
A Matlab multirate code that produces the results in Figure \ref{fig:wave}
is given in the Appendix.

\begin{figure}
\epsfig{file=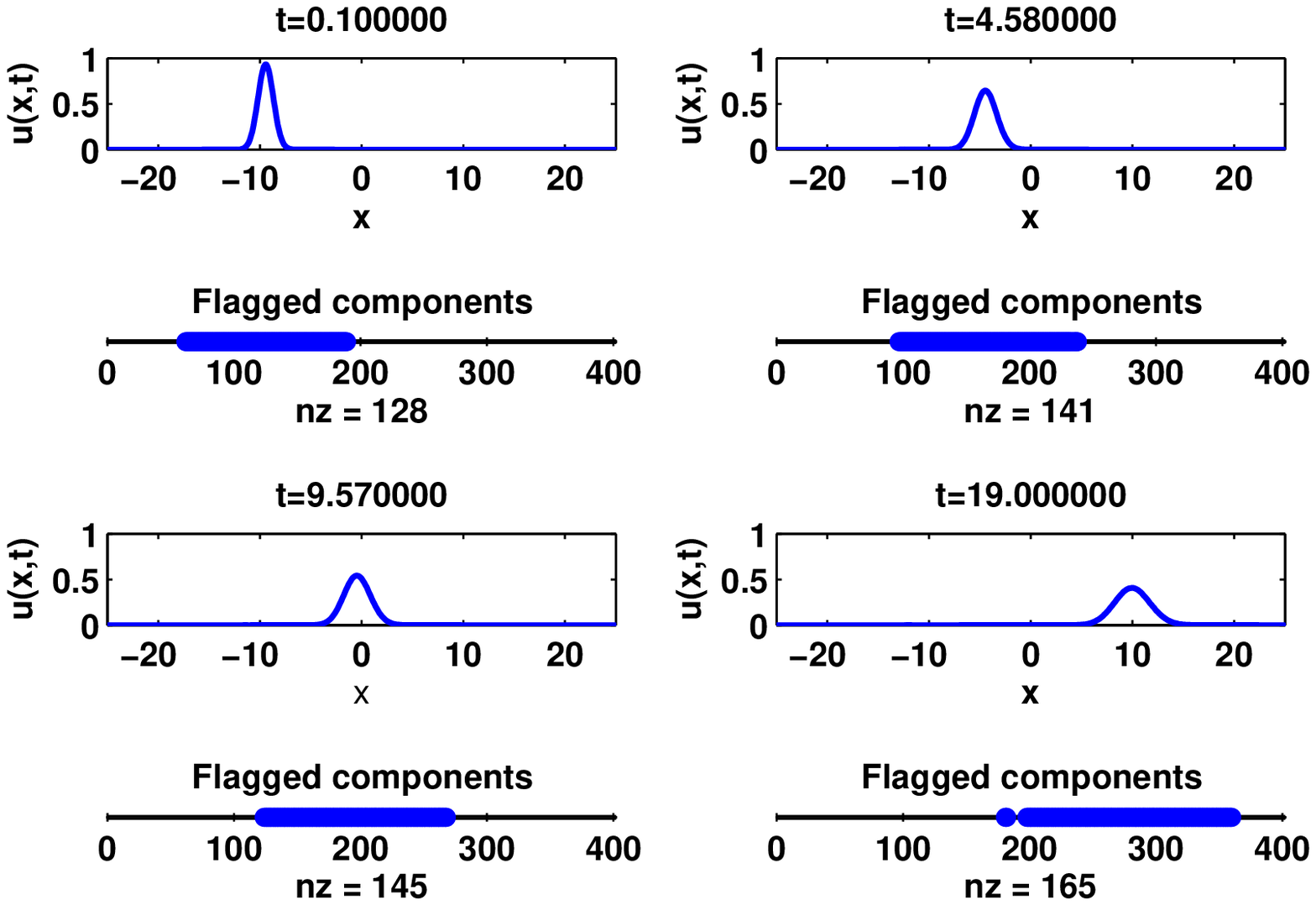,scale=0.8}
\caption{Solution of the advection equation, with unit wave speed,
at different times, obtained
using a multirate method. Components flagged as requiring a second integration
are shown in each case; the number of such components is 
given by the variable $nz$.
The total number of equations solved was $N=401$ and periodic
boundary conditions were used on the domain $-25 \leq x \leq 25$.
Multirate parameters (described in Algorithm \ref{algorithm2})
were $k=-6$ and $P=30$.}
\label{fig:wave}
\end{figure}

\subsection{Results}
\label{sec06:results}
Figure \ref{fig5}(a) shows the results of an integration with 
$\eps = 10^{-3}$. Note that only those steps which are near the facet tend to pack closely together,
but steps which are far away move relatively slowly and do not deviate significantly
from their initial uniform configuration. This expanding front of
closely packed steps represents the $t^{1/4}$ expansion of the 
facet radius \cite{IsraeliKandelPRB99,MargetisAzizPRB05}.

Figure \ref{fig5}(b) illustrates the separation in time scales
of the solution components and
shows which components of the solution are integrated
for a second time. 
As expected, our algorithm takes large time steps for components which are far
away from the facet. Near the facet and the collapsing top step, many relatively
small steps are taken. For the rapidly varying components
in this figure, only a representative sample of the meshpoints $t_n$ from the
$I_2$ integration are shown.

In contrast, when $\eps = 10^{-5}, 10^{-6}$, steps can be closely
packed even away from the facet. The plots in 
Figure \ref{fig6} show that a step bunching instability
arises when $\eps$ is sufficiently small and are qualitatively
very different to those in Figure \ref{fig5}.
The instability originates from steps with smaller radii 
and gradually spreads outwards so that more
and more steps bunch up. Our multirate scheme performs a second integration
when bunching and local stiffness arise: therefore,
our algorithm gradually becomes less efficient over time. However, as long as the
fraction of bunched steps is not too large, our algorithm remains competitive
compared to a standard adaptive 4th/5th order Runge-Kutta code. When
the fraction of bunched steps becomes close to unity, the optimal strategy
is to have the algorithm detect this automatically, and then switch to a fully
implicit, single-rate stiff solver. We leave this as future work, noting that 
inversion of the pentadiagonal Jacobian 
only costs $O(N)$ operations (where $N$
is the total number of existing steps).

\RRRfigOLD{
\EPfg{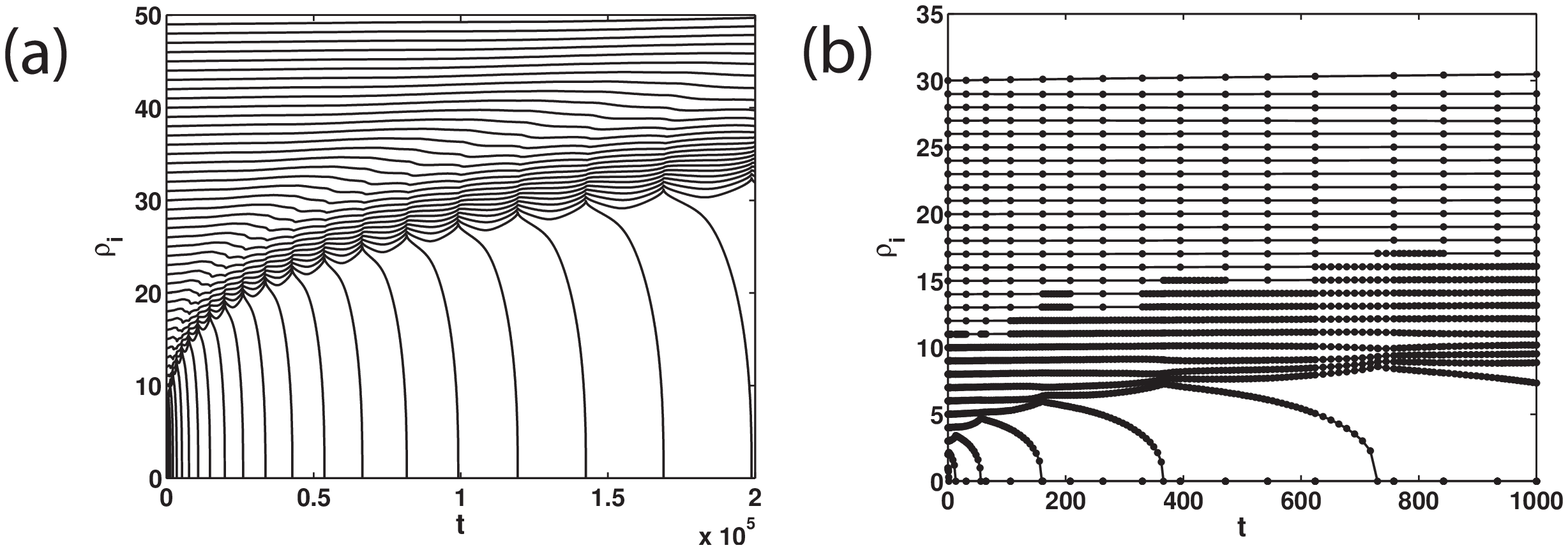,width=5.2in,height=2.2in} \hfill 
}
{\caption{A step instability arises when $\eps$ is sufficiently small. 
(a) Simulated relaxation of a nanostructure consisting of 200
steps, for a step train that is initially uniformly spaced. Only
the first 50 steps are shown. (b) Results of a simulation
with some of the time points of the multirate integration shown explicitly.}
\label{fig5}}

\RRRfigOLD{
\EPfg{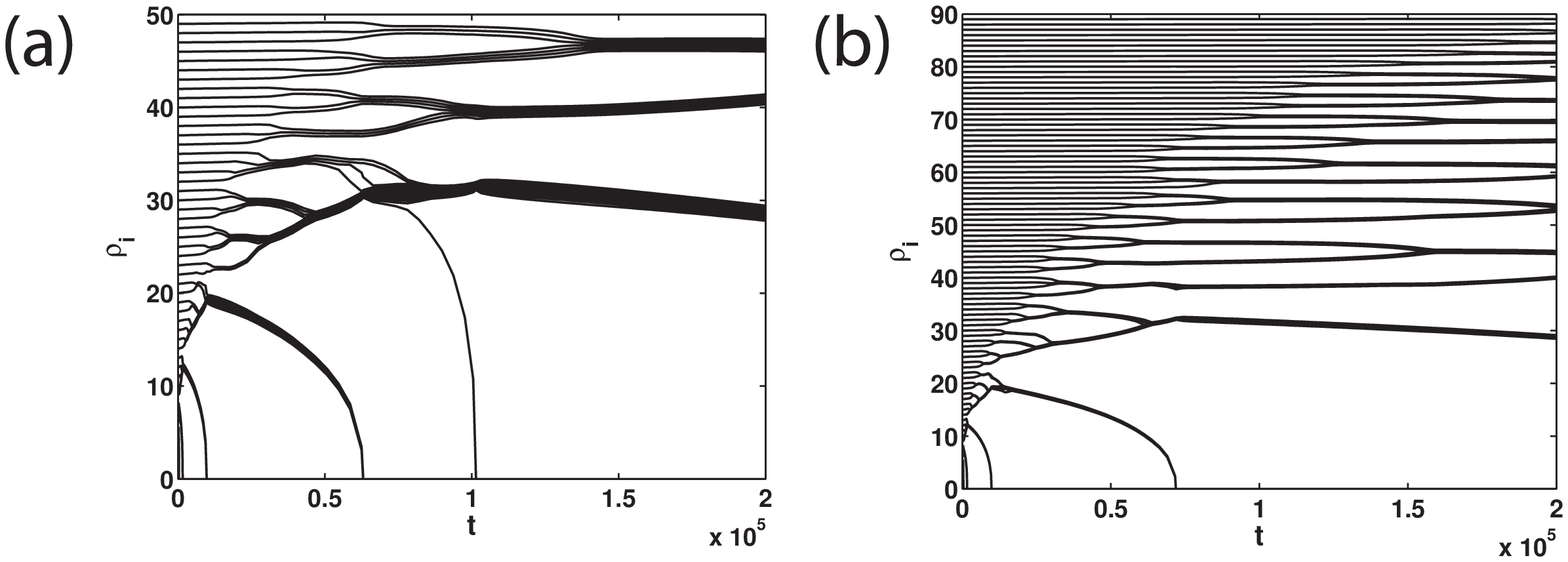,width=5.2in,height=2.1in} \hfill 
}
{\caption{A step bunching instability arises when $\eps$ is sufficiently small. 
(a) $\eps = 10^{-5}$.
(b) $\eps = 10^{-6}$. Approximately the first 50 and 90 steps
are shown respectively, out of 200. }
\label{fig6}}

\section{Conclusions}
\label{sec61:discussion}
In this paper, we present a multirate integration scheme 
that is designed to
efficiently solve the systems of ODEs that arise in the relaxation
of crystal mounds. These ODEs have two properties that
call for a multirate strategy: the singular collapse of
the innermost step and local stiffness. 
Our method
automatically detects singular/stiff components in the solution
and disregards them when computing
the size of the bulk (macro) time step. The result is that the
bulk timestep can be much larger than in a single rate method.
The trade-off is a re-integration of the stiff components
which usually consists of a small fraction of the total
number of equations in the ODE system.

Our method is globally fourth order when applied to
ODEs which have sufficiently smooth solutions -- for
example, the step equations studied 
in~\cite{SatoUwahaSurfSci99} and the wave equation
discretized through the method of lines in Section 
\ref{subsec:error-spread}. However,
the time singularities present in the axisymmetric step-flow
equations mean that near times of collapse, the integration of
steps near the facet suffers a loss in accuracy. 
To specifically deal with the singular inner trajectories,
our method couples a Simple Euler routine to the
bulk solver. Given that the truncation
error reduces to $O(\Delta t^{1/2})$ near the collapse
time, independent of the method order, Simple Euler
is the preferred method because it is computationally cheaper. 
Furthermore, the use of linear interpolation
to extract collapse times is consistent with the method's order.

The high order accuracy of our algorithm (when applied to bulk steps) 
relies on the ability to generate
high order interpolants during the run-time of a one-step integration
method. Our algorithm generates these interpolants by
using the intermediate stage function evaluations
of an embedded Runge-Kutta (RK) formula. Specifically, our
method computes $3^{rd}$ order
interpolants that are consistent with the $4^{th}$ order accuracy
of the integrator.
However, for general $n^{th}$ order RK formulae, we do not
know if it is always possible to construct interpolants
that have order $(n-1)$.

We see four main possible extensions to this work. 
The first is to generalize our multirate paradigm so that
it can be used for (i) higher order Runge-Kutta formulae
and (ii) multistep methods (e.g. BDF, Adams etc.)
We believe that it should be possible to make
\emph{any} method multirate -- the main obstacle
in doing this is to derive interpolants of a suitably
high order. 

The second is to explore in more detail the types of PDEs
that our multirate method can apply to. 
Typically, large systems of ODEs result from discretizations
of PDEs and it is for large ODE systems that our method
becomes competitive with single rate methods.
We think that
a basic requirement of the discretization is that it should
be strongly local. However, we have not fully explored which
discretizations are strongly local and which are not. For
example, a one-sided, upwind discretization 
$\dot{u}_n = (a/\Delta x) ( u_{n} - u_{n-1} )$
of the advection equation $u_t + a u_x = 0$ is strongly local
only when $a>0$. For $a<0$, a kronecker delta initial
condition is unstable and does \emph{not} remain localized.
A discretization using centered differences
$\dot{u}_n = (a/(2 \Delta x)) ( u_{n+1} - u_{n-1} )$
yields a system of ODEs that is never strongly local for any $a$.
For nonlinear equations, our method seems to be 
efficient for step-flow like ODEs with repulsive dipolar
step-step interactions. We were able to show that
the linearized step flow ODEs are strongly local; however,
we do not know if linearizing an ODE system is always 
sufficient to show strong locality.

The third is to explore how the choice of parameters
$k$ and $P$ in Algorithm \ref{algorithm2} affect the efficiency
of the integration and if there are optimal values of $k$ and $P$.
A ``good'' choice for $k$ and $P$ will
result in small number of re-integrated components and
a large macro-time step. If $k$ is too large 
and $P$ too small, the Algorithm behaves
like a single rate method. On the other hand, if $k$ is too small 
and $P$ is too large, many non-stiff components will be re-integrated along
with the stiff ones, rendering the method inefficient. 
Furthermore, our choice of $\mu$
as the $P^{th}$ percentile of the LTEs is somewhat arbitrary (but seems
to generate reasonable results). Another possibility
is to take the mean -- this amounts to increasing the sensitivity
of the bulk step size to the 
presence of one or two extremely stiff components. Clearly, the
performance of our multirate method is tied to the distribution
of LTEs, its moments, and identification of the ``largest 
errors''. Quantification of the
``largest errors'' and deciding
which moments to use is work in progress.
%

In summary, this work contributes to the currently
growing body of research in multirate methods.
We hope that
the strategies adopted in this paper can be
carried over to other physical problems and used
to improve the efficiency and accuracy of future multirate
algorithms.

\textbf{Acknowledgements} \\
We thank Dionisios Margetis for many helpful 
discussions and meetings. RRR was partially
supported by NSF grant DMS-0813648.

\section{Appendix}
Here we give the details of a Matlab multirate code to solve a wave equation.
\begin{verbatim}
function multirateCK

% solve the advection equation u_t + a u_x using a multirate method

global h N a

N = 401;
a = 1;
L = 50;

x = linspace(-L/2,L/2,N);
h = x(2)-x(1);
u0 = exp(-(x+10).^2)';
u = u0';

% parameters for integrator
desired_error = 1e-6; % desired error per step
MR = 1; % set MR=1 for multirate mode, MR=0 for single rate
T= 20;
t = 0;
dt = 0.1;
dt_max = 1.0;
k = -3;
P = 30; % approx percentage of components to reintegrate
safety = 0.95;

if MR == 1
    W = 10; % reintegrate more components on either side to be safe
else
    W = 0;
end

flags = zeros(1,N);
trynum=0;
numsteps=0;
r=round(N/2)+1;s=round(N/2)-1;

t1 = 0.1; t2 = 4.58; t3=9.57; t4=19;
plot_number=1;

while t<T
    [unew,error,K] = rk_onestep(t,u,dt,[1 N]);
    
    %%%%% start of multirate modification %%%
    if MR == 1
    
        % flag large errors
        flags = error > 10^k*percentile(error,P);        
        for i=1:length(flags)
            if flags(i) == 1
                r = i;
                break;
            end
        end
        for j=length(flags):-1:1
            if flags(j) == 1
                s = j;
                break;
            end
        end
        % just to be safe
        r = r-W;
        s = s+W;
        [unew,num_micro_steps] = micro_integrate(t,t+dt,u,unew,...
						                            [r s],K,desired_error*1e-3);
    end
    %%% end of multirate modification %%% 
    
    max_error = max( max(error(1:r-1)), max(error(s+1:end)) );
    R = ( desired_error/max_error )^(1/5);
    
    if R<1 % step failed
        dt = dt*safety*max(0.1,R); 
        trynum = trynum+1;
        if trynum>10
            sprintf('10 failed attempts!')
            return
        end
    elseif R>1 % step succeeded
        [r s (s-r) (s-r)/N]
        t = t+dt;
        u = unew;
        trynum = 0;
        numsteps = numsteps+1;     
        dt = dt*safety*min(5,R);
                
        if t+dt>T
            dt = (T-t);
        end
        
        if dt>dt_max
            dt = dt_max;
        end       
    end
    

    if t>t1 && plot_number == 1        
        subplot(4,2,1)
        plot(x,unew,'LineWidth',2);
        xlabel('x'); ylabel('u(x,t)');
        tit = sprintf('t=%f',t1); 
        title(tit);
        axis([-L/2 L/2 0 1]);
        subplot(4,2,3)
        spy(flags,20)
        title('Flagged components');
        plot_number = 2;
    elseif t>t2 && plot_number == 2
        subplot(4,2,2)
        plot(x,unew,'LineWidth',2);
        xlabel('x'); ylabel('u(x,t)');
        tit = sprintf('t=%f',t2); 
        title(tit);
        axis([-L/2 L/2 0 1]);
        subplot(4,2,4)
        spy(flags,20)
        title('Flagged components');
        plot_number=3;
    elseif t>t3 && plot_number == 3
        subplot(4,2,5)
        plot(x,unew,'LineWidth',2);
        xlabel('x'); ylabel('u(x,t)');
        tit = sprintf('t=%f',t3); 
        title(tit);
        axis([-L/2 L/2 0 1]);
        subplot(4,2,7)
        spy(flags,20)
        title('Flagged components');
        plot_number=4;
    elseif t>t4 && plot_number == 4
        subplot(4,2,6)
        plot(x,unew,'LineWidth',2);
        xlabel('x'); ylabel('u(x,t)');
        tit = sprintf('t=%f',t4); 
        title(tit);
        axis([-L/2 L/2 0 1]);
        subplot(4,2,8)
        spy(flags,20)
        title('Flagged components')
    end
    
end


%%%%%%%%%%%%%%%%%%%%%%%%%%%%%%%%%%%%%%%%%%%%%%%%%%%%%%%%%%%%%%%%%%%%%%%%%

function [y2,numsteps] = micro_integrate(t1,t2,y1,y2,cpt,K,desired_error)

% y(r) ... y(s) require integration
% y(r-2), y(r-1), y(s+1), y(s+2) are bcs.
% y2 requires updating

r = cpt(1); s = cpt(2);
dt = (t2-t1)/50;

y_current = y1;
t = t1;
trynum = 0;
numsteps = 0;
while t<t2
    
    [y_new,error] = rk_onestep(t,y_current,dt,[r-2,s+2]);     
    R = ( desired_error/max(error(r:s)) )^(1/5);
    if R<1 % step failed
        dt = dt*max(0.1,R); trynum = trynum+1;
        if trynum>10
            sprintf('10 failed attempts!')
            return
        end
    elseif R>1 % step succeeded
        t = t+dt; numsteps = numsteps+1;
        trynum = 0;
        
        y_new(r-2) = Interpolate(t,y1(r-2),y2(r-2),t1,t2,K(r-2,:));
        y_new(r-1) = Interpolate(t,y1(r-1),y2(r-1),t1,t2,K(r-1,:));
        y_new(s+1) = Interpolate(t,y1(s+1),y2(s+1),t1,t2,K(s+1,:));
        y_new(s+2) = Interpolate(t,y1(s+2),y2(s+2),t1,t2,K(s+2,:));
                
        y_current = y_new;
        dt = dt*min(5,R);          
        if t+dt>t2
            dt = (t2-t);
        end
        
    end
end
y2 = y_current;



%%%%%%%%%%%%%%%%%%%%%%%%%%%%%%%%%%%%%%%%%%%%%%%%%%%%%%%%%%%%%%%%%%%%%%%%%%%

function [ynew,error,K] = rk_onestep(t,y,dt,cpt)
% take a single rk45 with components cpt and 
% step of size dt and output the error

global h N a

a2 = 1/5; a3 = 3/10; a4 = 3/5; a5 = 1; a6 = 7/8;
b21 = 1/5;
b31 = 3/40; b32 = 9/40;
b41 = 3/10; b42 = -9/10; b43 = 6/5;
b51 = -11/54; b52 = 5/2; b53 = -70/27; b54 = 35/27;
b61 = 1631/55296; b62 = 175/512; b63 = 575/13824; 
b64 = 44275/110592; b65 = 253/4096;

c1 = 37/378;
c2 = 0;
c3 = 250/621;
c4 = 125/594;
c5 = 0;
c6 = 512/1771;

c1s = 2825/27648;
c2s = 0;
c3s = 18575/48384;
c4s = 13525/55296;
c5s = 277/14336;
c6s = 1/4;


k1 = dt*f(t,y,cpt);
k2 = dt*f(t+a2*dt,y+b21*k1,cpt);
k3 = dt*f(t+a3*dt,y+b31*k1+b32*k2,cpt);
k4 = dt*f(t+a4*dt,y+b41*k1+b42*k2+b43*k3,cpt);
k5 = dt*f(t+a5*dt,y+b51*k1+b52*k2+b53*k3+b54*k4,cpt);
k6 = dt*f(t+a6*dt,y+b61*k1+b62*k2+b63*k3+b64*k4+b65*k5,cpt);
K = [k1' k2' k3' k4' k5' k6'];

ynew_p = y + c1*k1 + c2*k2 + c3*k3 + c4*k4 + c5*k5 + c6*k6; % 5th order
ynew = y + c1s*k1 + c2s*k2 + c3s*k3 + c4s*k4 + c5s*k5 + c6s*k6; % 4th order

error = abs(ynew - ynew_p);
%%%%%%%%%%%%%%%%%%%%%%%%%%%%%%%%%%%%%%%%%%%%%%%%%%%%%%%%%%%%%%%

function ydot = f(t,y,components)
% evaluates RHS of ODE
% components = [r s] where 1 <= a1 < a2 <= N
% N = total number of ODEs

global h N a

r = components(1); s = components(2);

ydot = zeros(1,N);
ydot(r+1:s) = -(a/h)*(y(r+1:s) - y(r:s-1));


%%%%%%%%%%%%%%%%%%%%%%%%%%%%%%%%%%%%%%%%%%%%%%%%%%%%%%%%%%%%%%%%%

function yinterp = Interpolate(t,y1,y2,t1,t2,k)
a = (t-t1)/(t2-t1);
yinterp = y1+a*k(1)+0.5*a.^2*(-8/3 * k(1)+25/6*k(4)-3/2*k(5))+ ...
            a.^3/6*(10/3*k(1) - 25/3 * k(4) + 5*k(5));


function out = percentile(X,n)
% outputs the nth percentile for data X.
% e.g. n = 50 ---> out = median
% e.g. n = 25 ---> out = X* such that 75% of X are smaller than X*

X = sort(X,'descend'); % largest to smallest
J = round(0.01*n*length(X));
out = X(J);

\end{verbatim}

%

\end{document}